\documentclass[final,1p,times]{elsarticle}

\usepackage[utf8]{inputenc}
\usepackage[T1]{fontenc}
\usepackage{makecell} % long text in table cells
\usepackage{booktabs} % for prettier tables
\usepackage{array}    % for table column formatting
\usepackage{subcaption}
\usepackage{algpseudocode}
\usepackage{algorithm}
\usepackage{fancyhdr}
\usepackage{comment}
\usepackage{tikz,graphicx}
\usepackage{multirow}
\usepackage{amsfonts,amsmath,amssymb}
\usepackage{mathtools}
\usepackage{mathrsfs}
\usepackage{upgreek}

\newcommand{\sR}{{\mathbb{R}}}

\newcommand{\bx}{{\mathbf{x}}}

\newcommand{\bz}{{\mathbf{z}}}
\newcommand{\bu}{{\mathbf{u}}}
\newcommand{\bv}{{\mathbf{v}}}

\newcommand{\bo}{{\mathbf{0}}}

\newcommand{\bI}{{\mathbf{I}}}

\newcommand{\rv}[1]{\boldsymbol{#1}}

\newcommand{\bn}{{\bf n}}
\newcommand{\sumin}{\sum_{K \in \mathcal{K}_h}}

\begin{document}

\begin{frontmatter}

\title{Ensemble Score Filter for Data Assimilation of Two-Phase Flow Models in Porous Media}

\author[math]{Ruoyu Hu}
\author[sc]{Sanjeeb Poudel}
\author[math]{Feng Bao}
\author[math]{Sanghyun Lee}

\affiliation[math]{organization={Department of Mathematics, Florida State University},%
            city={Tallahassee},
            postcode={32306}, 
            state={Florida},
            country={USA}}
\affiliation[sc]{organization={Department of Scientific Computing, , Florida State University},
            city={Tallahassee},
            postcode={32306}, 
            state={Florida},
            country={USA}}
%\ead[cor]{slee17@fsu.edu}
\begin{abstract}
Numerical modeling and simulation of two-phase flow in porous media is challenging due to the uncertainties in key parameters, such as permeability.  To address these challenges, we propose a computational framework by utilizing the novel  Ensemble Score Filter (EnSF) to enhance the accuracy of state estimation for two-phase flow systems in porous media. The forward simulation of the two-phase flow model is implemented using a mixed finite element method, which ensures accurate approximation of the pressure, the velocity, and the saturation. The EnSF leverages score-based diffusion models to approximate filtering distributions efficiently, avoiding the computational expense of neural network-based methods. By incorporating a closed-form score approximation and an analytical update mechanism, the EnSF overcomes degeneracy issues and handles high-dimensional nonlinear filtering with minimal computational overhead.  Numerical experiments demonstrate the capabilities of EnSF in scenarios with uncertain permeability and incomplete observational data. 
%  This approach offers a robust, efficient solution for advancing data assimilation in porous media applications, \textcolor{cyan}{including reservoir management and environmental risk assessment.} \\ \textcolor{cyan}{We do not show any of these applications for our approach. }
\end{abstract}

\end{frontmatter}

\section{Introduction}

%\cite{evensen1994sequential}

%\textbf{Two Phase Flow }

The study of porous media holds significant importance given its wide-ranging applications from oil and gas exploration to geo-energy systems, including geothermal energy, CO$_2$ sequestration, and the transport of contaminants in the groundwater system. 
The numerical modeling and simulation techniques of porous media have gained traction and proven invaluable for understanding and predicting subsurface behavior~\cite{Jim1990,Mary}. However, uncertainties in critical parameters, such as permeability tensor, often compromise the solution’s reliability. In addition, the stochastic nature of geological formations exacerbates these challenges, leading to substantial variability in model outputs~\cite{watanabe2010uncertainty,zhang2006stochastic}.

This study explores the application of recently developed novel data assimilation techniques based on Ensemble Score Filter (EnSF)
to improve the accuracy of numerical solutions for the two immiscible fluids in porous media. Data assimilation plays a crucial role in bridging the gap between model predictions and observation \cite{law2015data}. By systematically integrating observed data, data assimilation techniques enable the refinement of state variables estimates, while accounting for uncertainties in both the model parameters and the observational data. This is particularly vital for two-phase flow systems, where the state of the system evolves nonlinearly and is highly sensitive to spatial variability in model parameters. Effective data assimilation not only improves the predictive accuracy of the model but also enhances its utility in decision-making processes for critical applications such as reservoir management, geothermal energy extraction, and environmental risk assessment \cite{carrassi2018data,rubin1991transport}. 

%In the literature, several data assimilation techniques have been proposed to tackle the challenges associated with model uncertainties and observational limitations. Among these, continuous data assimilation has emerged as a promising approach~\cite{azouani2014continuous}, particularly for modeling two-phase flow in porous media~\cite{chow2022continuous}, due to its ability to systematically incorporate observational data into the model framework over time. Continuous data assimilation has been widely adopted in theoretical studies because of its capability to provide rigorous analytical results and stability guarantees under appropriate conditions. However, its practical application often encounters challenges arising from the complex dynamics and heterogeneity inherent in subsurface environments. \textcolor{cyan}{Some limitations of continuous data assimilation that can be overcome using EnSF}
%{\color{red}{The key feature of the proposed algorithm in this paper, which differentiates it from the current state of the art, is its emphasis on addressing the fact that data may often be incorrect or that errors in the data may go unrecognized, leading to mismatches with real-world observations.}}

The primary mathematical tool for addressing data assimilation problems is \textit{optimal filtering}. 
The goal of the optimal filtering problem is to approximate the conditional probability density function (PDF) of the state process given observational data, referred to as the \textit{filtering distribution} ~\cite{Bao_Zakaid_2015, vanLeeuwen_2009}. The corresponding conditional expectation provides the optimal state estimation. Bayesian inference plays a central role in solving optimal filtering problems, with Kalman-type filters being a prominent example. For linear state dynamics and observations under a Gaussian assumption, the Kalman filter provides an analytical solution.

To handle nonlinear systems, Ensemble Kalman filters (EnKF) were developed \cite{evensen1994sequential,houtekamer1998data}. However, EnKF assumes Gaussian filtering densities, which may not hold for highly nonlinear systems or observations. Beyond EnKF, several effective methods have been proposed for nonlinear optimal filtering, including the particle filter \cite{ andrieu2010particle,chorin2009implicit,gordon1993novel,kang2018improved,pitt1999filtering,snyder2008obstacles,bao2014hybrid, Sumner2025Nano}, Zakai filter \cite{Bao_Zakaid_2015, zakai1969optimal}, and methods based on stochastic partial differential equations (SPDEs) \cite{bao2021data,bao2015forward,bao2016first,bao2017adaptive,hu2002approximation,gobet2006discretization}. The particle filter, also known as the sequential Monte Carlo method, is well-suited for capturing complex non-Gaussian filtering densities in moderate dimensions. However, its performance degrades in high-dimensional problems due to the challenges of high-dimensional sampling in Bayesian inference. SPDE-based methods, such like Zakai type filters \cite{zakai}, while mathematically rigorous and capable of providing stable state estimation over time, face significant computational challenges. Solving SPDEs becomes increasingly expensive as the problem's dimensionality grows, with computational costs escalating exponentially \cite{Bao_Zakaid_2015, zhang2008grid,cai2020learning,dhariwal2021diffusion}. In particular, in the case of two-phase flow, the coupled variables of velocity, pressure, and saturation lead to a high dimensional data assimilation problem.

%\textcolor{cyan}{
%Continuous data assimilation~\cite{azouani2014continuous} has been studied for modeling two-phase flow in porous media~\cite{chow2022continuous}, due to its ability to systematically incorporate observational data into the model framework over time. Additionally, it has been widely adopted in theoretical studies because of its capability to provide rigorous analytical results and stability guarantees under appropriate conditions. However, its practical application often encounters challenges arising from the complex dynamics and heterogeneity inherent in subsurface environments.}
%{\color{blue} Novelties B: EnSF and the applicability of EnSF in practical scenarios.}

In this paper, we address the challenge of state estimation for solutions of the two-phase flow model by introducing a novel Ensemble Score Filter (EnSF) designed for accurate and efficient data assimilation. In the literature, continuous data assimilation~\cite{azouani2014continuous} has been studied for modeling two-phase flow in porous media~\cite{chow2022continuous}. However, continuous data assimilation operates under the perfect model assumption, which neglects intrinsic uncertainties in PDE models. This omission leads to inconsistencies between practical data and the target PDE models, as real-world systems often deviate from idealized mathematical formulations. The EnSF relies on score-based diffusion models \cite{NEURIPS2021_49ad23d1, NEURIPS2020_4c5bcfec, NEURIPS2019_3001ef25, song2021scorebased} to model the evolution of filtering distributions in a pseudo-temporal domain. Traditional diffusion-based filtering methods \cite{bao2023scorebased} rely on neural networks to estimate score functions, requiring time-intensive retraining at each filtering step and significant storage for high-dimensional problems. EnSF addresses these challenges by entirely avoiding neural network training. Instead, it directly approximates the score function using a closed-form expression combined with mini-batch Monte Carlo estimators, enabling efficient and accurate score estimation. Furthermore, EnSF incorporates an analytical update step to progressively integrate observational data, reducing degeneracy issues in high-dimensional nonlinear filtering. Extensive studies demonstrate that EnSF outperforms state-of-the-art methods \cite{Toan2025online, liang2024online, EnSF_2023, Bao_SC_2024, Inpainting_2024}, such as the Local Ensemble Transform Kalman Filter (LETKF), in various scenarios, including highly challenging high-dimensional data assimilation problems for forecasting chaotic atmospheric models \cite{Bao_MWR_2024}.

%To demonstrate the effectiveness of the EnSF in data assimilation for state estimation in the two-phase flow model,  we design experiments that closely replicate practical scenarios. Specifically, we assume the model incorporates uncertainties in permeability arising from incomplete knowledge, with no explicit mathematical representation of these uncertainties due to the difficulty of accurately determining permeability over large regions. Moreover, we assume that observational data for the two-phase flow model are obtained through nonlinear detectors and are spatially incomplete across the domain.

To demonstrate the effectiveness of the EnSF in data assimilation for state estimation in the two-phase flow model,  we design experiments that closely replicate practical scenarios. Specifically, we assume the model incorporates uncertainties in permeability arising from incomplete knowledge. Additionally, we assume that the observational data is sparse, nonlinear, and noisy. In the experiments, we show that EnSF for data assimilation performs well in high-dimensional problems highlighting its efficient computational complexity. Additionally, we compare the results with LETKF, which also demonstrates that our method outperforms LETKF.

The remainder of this paper is organized as follows. In Section \ref{sec:2}, we present the mathematical model for two-phase flow in porous media, along with a detailed discussion of the general approach to data assimilation. In Section \ref{sec:3}, we present the details of the novel score-based diffusion model, highlighting its theoretical underpinnings and potential advantages for capturing the complex dynamics of multi-phase flow systems. This section also explores the innovative aspects of the model that distinguish it from existing methodologies.
Section \ref{sec:4} introduces the overall coupled computational framework, which integrates the finite element solver with the Ensemble Smoother with Sampling (EnSF) algorithm. This section provides an in-depth explanation of the framework’s architecture, its implementation, and its role in achieving reliable and efficient numerical solutions for data assimilation problems.
Finally, in Section \ref{sec:5}, we present the results of numerical experiments designed to illustrate the capabilities and performance of the proposed algorithms. Through these experiments, we demonstrate the effectiveness of our approach in addressing the challenges of modeling two-phase flow in porous media and highlight its practical applicability across various scenarios.

\section{Data assimilation for the two-phase flow model}
\label{sec:2}

In this section, we present the data assimilation model for the two-phase flow. We begin with a brief overview of the two-phase flow model in porous media. Next, we introduce the data assimilation framework, which aims to estimate the state of a hidden dynamical system using partial and noisy observational data.

\subsection{The two-phase flow model in porous media}
%\vspace{0.5em}

%The fluid velocity in porous media is described by Darcy's law \cite{whitaker_flow_1986}. According to this law, the velocity is proportional to the pressure gradient. Similarly, 
For immiscible and incompressible two-phase flow in porous media, the fluid velocity is described based on Darcy's law \cite{whitaker_flow_1986, chow_continuous_2022, lee_enriched_2018, li2009upscaling}.
Let $\Omega \in \mathbb{R}^d$ be a bounded domain for $d = 2, 3$ with simply connected Lipschitz boundary $\partial \Omega$, and let $T>0$ be the given final time.  The velocity of phase $j$ $(j = 1, 2)$, is denoted by $\tilde{\bu}_j : \Omega \times (0, T] \rightarrow \mathbb{R}^d$ and is given by
\begin{equation} 
    \tilde{\bu}_j := -\lambda_j(s) \mathbf{K}\nabla p_j, \quad \text{in } \Omega \times (0, T],    
\end{equation}
where $s : \Omega \times (0, T] \rightarrow \mathbb{R}$ is the scalar-valued saturation, and $p : \Omega \times (0, T] \rightarrow \mathbb{R}$ is the scalar-valued pressure.
Here, $\lambda_j(s) := {k_{rj}(s)}/{\mu_j}$ is the mobility, with $k_{rj}(s)$ being the relative permeability as a function of saturation, 
%(which represents the volume fraction of the phase and ranges from 0 to 1),
$\mu_j$ is the viscosity of phase $j$, and ${\bf K} = K {\bf{I}} \in \mathbb{R}^{d\times d}$ is the absolute permeability tensor. $K$ is a constant and 
 ${\bf{I}}$ is the identity matrix. For simplicity, we assume zero capillary pressure, i.e., $p_1 = p_2$, and denote the pressure with $p$. 
Thus, from the conservation of mass for each phase, we obtain
\begin{equation}
    \nabla \cdot \tilde{\bu}_j = q_j, \quad \text{in } \Omega \times (0, T], 
\end{equation}
where $q_j$ represents the volumetric source term for each phase. After summing the two phases, we get the pressure equation, 
\begin{equation}
    -\nabla \cdot ({\bf K} \lambda(s) \nabla p) = q,
\end{equation}
where $q = q_1 + q_2$ is the total volumetric source term.
Here, $\lambda (s)$ is the total mobility that is given by
\begin{equation}
    \lambda(s) = \lambda_1(s) + \lambda_2(s).
\end{equation}
For a constant porosity $\phi$, the saturation of phase $j$, denoted by $s_j$, is governed by the conservation equation
\begin{equation}
    \phi \frac{\partial s_j}{\partial t} + \nabla \cdot \tilde{\bu}_j = q_j, \quad \text{in } \Omega \times (0, T].
\end{equation}
For the phase $j = 1$, we define the fraction flow $F(s)$ as
\begin{equation}
    F(s) := \frac{\lambda_1 (s)}{\lambda_1(s) + \lambda_2 (s)}.
\end{equation}
With the above definition of the fractional flow~\cite{pope1980application}, the saturation of phase $j=1$ can be written as
\begin{equation}
    \phi \frac{\partial s_1}{\partial t} + \nabla \cdot \left( F(s) \tilde{\bu} \right) = q_1, \quad \text{in } \Omega \times (0, T].
\end{equation}
With $q_1 = F(s) q$, $\bu = \tilde{\bu}/ \phi$, and denoting $s_1$ with $s$, the saturation equation for phase $j = 1$ can be written as 
\begin{equation}
    \frac{\partial s}{\partial t} +
    \bu \cdot \nabla F(s) = {0}, \quad \text{in } \Omega \times (0, T].
\end{equation}
In summary, we have the following two equations for the pressure and saturation as our governing system:
\begin{subequations}
\label{eqn:main}
    \begin{alignat}{2}
        -\nabla \cdot ({\bf K} \lambda(s) \nabla p) &= q, \quad \text{in } \Omega \times (0, T], \\
        \frac{\partial s}{\partial t} +
        \bu \cdot \nabla F(s) &= {0}, \quad \text{in } \Omega \times (0, T],
    \end{alignat}
\end{subequations}
where the above system is supplemented by the following boundary (using Dirichlet conditions) and initial conditions:
\begin{subequations}
    \begin{alignat}{2}
        p &= p_D, \quad &&\text{on } \partial \Omega \times (0, T], \\
        s &= s_D, \quad &&\text{on } \partial \Omega \times (0, T], \\
        s(\cdot, 0) &= s^0, \quad &&\text{in } \Omega.
    \end{alignat}
\end{subequations}

In practical applications, these solvers face significant challenges arising from uncertainty. For instance, the permeability data, denoted as ${\bf{K}}$, may be unknown and unobservable in the two-phase flow model. This lack of information can introduce unknown and unquantifiable model errors due to incomplete physical knowledge. Consequently, these errors may cause the PDE system \eqref{eqn:main} to deviate from the true underlying physics, potentially undermining the reliability of the simulation results.

To address the model and simulation errors mentioned above, field-collected data can be utilized to calibrate and refine the numerical solution, thereby reducing discrepancies between the model predictions and the actual physical system. This process naturally leads to the formulation of a data assimilation problem, where observed data are systematically integrated into the computational framework to enhance model accuracy. In the following subsection, we present the general framework of data assimilation, outlining its principles and relevance to improving the reliability of two-phase flow simulations.

\subsection{The data assimilation problem}

In the data assimilation problem, we consider the following state dynamics
\begin{equation}\label{NLF:State}
X_{t_{n+1}} = f(X_{t_{n}}, \omega_{t_{n}}), \qquad n = 0, 1, 2, \cdots 
\hspace{2em} \text{(State)}
\end{equation}
where $X_{t_{n}}$ is the state of a quantity of interest at time instant $t_n$, 
$f(\cdot)$ is a forward dynamical model, which contains a random variable $\omega_{t_n}$. 
For instance, in our work, $X_{t_n}$ is the pressure, velocity or the saturation for the two-phase flow model, and the random variable $\omega_{t_n}$ is employed as an input for $K(\cdot, \omega_{t_n})$ and $s^0(\cdot, \omega_{t_n})$.

%The dynamical model in this work is the PDE system Eq. \eqref{eqn:main}, and the random variable $\omega_{t_n}$ represents the uncertainty in the model  
%{\color{red} which perturbs the PDE input data such as permeability or initial conditions.}
%The specific probability distribution of $\omega_{t_n}$ may depend on the source of uncertainty we consider in the model. It is important to note that, in practical implementations, the dynamical model $f(\cdot)$ is typically approximated using a numerical solver, and {\color{red} we utilize low-order locally mass conservative finite element methods in this paper.} \textcolor{cyan}{Will it be better to have this paragraph in section 4? We have used $\bar{f}$ for f.}

The goal of the data assimilation problem is to utilize partial noisy observational data on $X$ to reduce uncertainty in the state model $f$ and make the best estimate of $X$. We denote the observations as $Y_{t_{n+1}}$, given by
\begin{equation}\label{NLF:Obs}
Y_{t_{n+1}} = h(X_{t_{n+1}}) + \epsilon_{t_{n+1}}, \hspace{2em} \text{(Observation)}
\end{equation}
where $h(\cdot)$ is the observation function that may provide only \textit{partial and indirect measurements} of $X$, and $\epsilon_{t_{n+1}}$ represents observational noise.  

The process of finding the best estimate for $X_{t_{n+1}}$ is known as the 
\textit{optimal filtering problem}. 
The optimal filter, denoted by 
\begin{equation}
\tilde{X}_{t_{n+1}} := \mathbb{E}[X_{t_{n+1}}\big| Y_{1:t_{n+1}}],
\end{equation}
is the conditional expectation of $X_{t_{n+1}}$ given the $\sigma$-algebra $Y_{1:t_{n+1}}$, which contains all observational information up to time $t_{n+1}$. In practice, the conditional expectation is not directly approximated. 
Instead, the focus is typically on approximating the conditional probability density function (PDF) of the state, denoted by $P(X_{t_{n+1}}|Y_{1:t_{n+1}})$, commonly referred to as the filtering distribution.

Since the state model \eqref{NLF:State} involves uncertainty and the observation model \eqref{NLF:Obs} provides only partial and indirect information about the state $X$, relying solely on either the model or the data is insufficient to produce accurate estimates of the true state of the target system. The standard approach to solving the data assimilation problem is the Bayesian filter, which recursively applies Bayesian inference to dynamically integrate the model and data. The general framework of recursive Bayesian filter consists of two steps at each time interval 
$[t_n, t_{n+1}]$: a prediction step and an update step. 

In the prediction step, assuming that the filtering density $P(X_{t_n}|Y_{1:t_n})$ is obtained at the time instant $t_n$, we use the Chapman-Kolmogorov formula to propagate the state equation in \eqref{NLF:State} from $t_n$ to $t_{n+1}$ and obtain the prior filtering distribution as 
\begin{equation}\label{Prediction}
    \textit{Prior filtering distribution:} \quad P(X_{t_{n+1}}|Y_{1:t_n}) = \int P(X_{t_{n+1}}|X_{t_{n}})P(X_{t_{n}}|Y_{1:t_n})dX_{t_n}, 
\end{equation}
where $P(X_{t_{n+1}}|X_{t_{n}})$ is the transition probability derived from the state dynamics described by Eq. \eqref{NLF:State}. The prior filtering distribution $P(X_{t_{n+1}}|Y_{1:t_n})$ contains the data information up to time $t_{n}$ and the model information captured by the transition probability.

In the update step, we obtain the posterior filtering distribution by incorporating the new observational data $Y_{t+1}$ with the prior filtering distribution through the following Bayesian inference formula: 
\begin{equation}\label{Update}
    \textit{Posterior filtering distribution:} \quad P(X_{t_{n+1}}|Y_{1:t_{n+1}}) \propto P(X_{t_{n+1}}|Y_{1:t_n})P(Y_{t_{n+1}}|X_{t_{n+1}}), 
\end{equation}
where the likelihood function $P(Y_{t_{n+1}}|X_{t_{n+1}})$ is defined by 
\begin{equation}\label{Likelihood}
    P(Y_{t_{n+1}}|X_{t_{n+1}}) \propto \text{exp}\big[-\frac{1}{2}(h(X_{t_{n+1}})-Y_{t_{n+1}})^\top R^{-1}(h(X_{t_{n+1}})-Y_{t_{n+1}})\big], 
\end{equation}
with $R$ denoting the covariance matrix of the random noise $\epsilon$ in \eqref{NLF:Obs}.

In the following section, we introduce our novel ensemble score filter for solving the data assimilation problem within the diffusion model-based generative artificial intelligence (AI) framework.

\section{The training-free approach for the score-based diffusion model}
\label{sec:3}

The ensemble score filter adopts the diffusion model, which is one of the most important generative AI techniques, to generate samples that follow the filtering distributions desired in the data assimilation problem. In this section, we shall introduce a training-free approach for generative sampling in the score-based diffusion model. 

\subsection{The score-based diffusion model}

In a diffusion model, we have the following $\sR^{d}$-dimensional stochastic differential equation (SDE)
\begin{equation}\label{forward:SDE}
d \rv{Z}_{\tau} = b(\tau)\rv{Z}_{\tau} d\tau + \sigma(\tau)d\rv{W}_{\tau} \quad \text{(forward SDE),}
\end{equation}
where $\rv{W}_{\tau} \in \sR^{d}$ is a standard Brownian motion (Wiener process) corresponding to an It\^o integral $\int \cdot d\rv{W}_{\tau}$, while $b$ and $\sigma$ are two explicitly given functions referred to as the drift and diffusion coefficients, respectively. The initial condition $\rv{Z}_0$ of the SDE in Eq.~\eqref{forward:SDE} follows a target distribution with its PDF denoted by $Q_{0}(\bz_0)$. One can show that with properly chosen $b$ and $\sigma$, the diffusion process $\{\rv Z_{\tau}\}_{0\leq \tau \leq 1}$ can transform any target PDF ${Q}_0(\bz_0)$ to a standard Gaussian, i.e. $\rv Z_1 \sim \mathcal{N}(\bo, \bI_d)$, where $[0, 1]$ is often referred to as the \textit{pseudo-time interval}. While there are multiple options for $b$ and $\sigma$, in this work we let
\begin{equation}\label{coefficients}
\begin{aligned}
    &b(\tau) = \frac{d\log \alpha_{\tau}}{d\tau}, \qquad \sigma^2(\tau) = \frac{d \beta^2_{\tau}}{d \tau} - 2 \frac{d \log \alpha_{\tau}}{d \tau} \beta^2_{\tau}
\end{aligned}
\end{equation}
with $\alpha_{\tau} = 1 - \tau$ and $\beta_{\tau} = \sqrt{\tau}$ for $\tau \in [0, 1]$, which is consistent with the choice made in \cite{song2021scorebased}. 

Generating samples $\{\rv z^{(m)} \}_{m=1}^{M}$ of the target random variable $\rv Z_0$ pertains to simulating the following reverse-time SDE:
\begin{equation}\label{reverse:SDE}
\begin{aligned}
d\rv Z_{\tau} = \left[ b(\tau)\rv Z_{\tau} -\sigma^2(\tau) \mathcal{S}(\rv Z_{\tau},\tau) \right] dt &+ \sigma(\tau)d\overleftarrow{\rv W}_{\tau} \quad \text{(reverse-time SDE),}
\end{aligned}
\end{equation}
where $\int \cdot d\overleftarrow{\rv W}_{\tau}$ is a \textit{backward} It\^o stochastic integral \cite{SDE, BDSDE}, and $\mathcal{S}(\cdot, \tau)$ is the so-called \textit{score function} given by

\begin{equation}\label{eq:score_func}
\mathcal{S}(z_{\tau}, \tau) \coloneqq \nabla \log Q_{\tau}(\bz_{\tau}),
\end{equation}
where $Q_{\tau}(\bz_{\tau})$ denotes the PDF of $\rv Z_{\tau}$.
Note that the reverse-time SDE in Eq.~\eqref{reverse:SDE} is also a diffusion process except that the propagation direction is backwards in time from $1$ to $0$ with initial condition $\rv Z_1$ given at time $1$. An important result from the literature is that the solution $\rv Z_0$ of the reverse-time SDE follows the target distribution \cite{SF_2023}. The practical implications are that we can generate samples from a standard Gaussian distribution (which can be done efficiently) and use the reverse-time SDE to transform them to samples of the target distribution. The score function $\mathcal{S}(\cdot, \tau)$ has an important role in this mapping process as it stores information about the distribution of the samples, which in turn helps guide their transformation over the pseudo-time interval as $\tau \rightarrow 0$. In particular, having the score function associated with the target PDF $Q_0(\bz_0)$ and the predefined forward SDE allows us to generate an unlimited number of target samples by running the reverse-time SDE in Eq.~\eqref{reverse:SDE}.

The traditional use of diffusion models in the machine learning (ML) literature is to generate highly realistic images and videos \cite{song2021scorebased}. However, these image/video processing-oriented diffusion models require extensive neural network training in order to estimate the score function $\mathcal{S}(z_{\tau}, \tau)$. In data assimilation applications, the score function must be estimated repeatedly within the recursive Bayesian inference framework, which makes deep learning-based approximations impractical due to the high computational cost of neural network training. To address this challenge, we introduce a Monte Carlo-based ensemble approximation for the score function, which serves as the foundation for our ensemble score filter algorithm in data assimilation.

\subsection{Monte Carlo based approximation for the score function}
Note that the forward process Eq.\eqref{forward:SDE} is a linear SDE.
The definitions for the drift and diffusion coefficient in Eq.~\eqref{coefficients} can ensure that the conditional density function $Q_{\tau}(\rv Z_\tau | \rv Z_0)$ for any fixed $\rv Z_0$ is the following Gaussian distribution:
\begin{equation}\label{eq:gauss}
Q_{\tau}(\rv Z_\tau | \rv Z_0) = \mathcal{N}(\alpha_{\tau} \rv Z_0, \beta_{\tau}^2 \mathbf{I}_d), 
\end{equation}
Substituting $Q_\tau(\rv Z_\tau)= Q_\tau(\rv Z_\tau, \rv Z_0) = Q_\tau(\rv{Z}_\tau |\rv Z_0) Q_\tau(\rv Z_0)$ into Eq.~\eqref{eq:score_func} and exploiting the fact in Eq.~\eqref{eq:gauss}, we can rewrite the score function as 
\begin{equation}\label{eq:score11}
\begin{aligned}
\mathcal{S}(\rv Z_{\tau}, \tau) & = \nabla_z \log \left(\int_{\mathbb{R}^d} Q_\tau(\rv{Z}_\tau | \rv Z_0) Q_\tau(\rv Z_0) d\rv Z_0\right)\\
& = \frac{1}{\int_{\mathbb{R}^d} Q_\tau(\rv{Z}_\tau | \rv Z'_0) Q_\tau(
\rv Z'_0) d\rv Z'_0}   \int_{\mathbb{R}^d}  - \frac{\rv Z_\tau - \alpha_\tau \rv Z_0}{\beta^2_\tau}Q_\tau(\rv{Z}_\tau | \rv Z_0) Q_\tau(
\rv Z_0) d\rv Z_0\\
& =  \int_{\mathbb{R}^d}  - \frac{\rv Z_\tau- \alpha_\tau \rv Z_0}{\beta^2_\tau} w_\tau(\rv{Z}_\tau,  \rv Z_0)  Q_0(\rv Z_0)d\rv Z_0,\\
\end{aligned}
\end{equation}
where the weight function $w_\tau(\rv{Z}_\tau, \rv Z_0)$ is defined by
\begin{equation}\label{eq:weight}
w_\tau(\rv{Z}_\tau,  \rv Z_0) :=  \frac{ Q_\tau(\rv{Z}_\tau | \rv Z_0) }{\int_{\mathbb{R}^d} Q_\tau(\rv{Z}_\tau | \rv Z'_0) Q_\tau(\rv Z'_0) d\rv Z'_0},
\end{equation}
satisfying that $\int_{\mathbb{R}^d}w_\tau(\rv{Z}_\tau, \rv Z_0) d\rv Z_0 = 1$. 

Then, we introduce the following Monte Carlo approximation for the score 
\begin{equation}\label{eq:MC}
\mathcal{S}(z, \tau) \approx \bar{\mathcal{S}}(z, \tau) :=  \sum_{j=1}^{J} - \frac{z - \alpha_\tau z^{(j)}}{\beta^2_\tau} \bar{w}_\tau({z},  z^{(j)})), 
\end{equation}
using a mini-batch with batch size $J \le M$, denoted by $\{z^{(j)}\}_{j=1}^J$, from the ensemble of original data samples $\{z^{(m)}\}_{m=1}^M$,  and the weight $w_\tau({z},  z^{(j)})$ is calculated by
\begin{equation}\label{eq:weight_app}
w_\tau({z},   z^{(j)}) \approx  \bar{w}_\tau({z},   z^{(j)}) := \frac{Q_\tau(z |  z^{(j)}) }{\sum_{j=1}^{J} Q_\tau(z| z^{(j)})},
\end{equation}
This means $w_\tau({z},  z^{(j)})$ can be estimated by the normalized probability density values $\{Q_{\tau}(z|z^{(j)})\}_{j=1}^J$. In practice, the mini-batch $\{z^{(j)}\}_{j=1}^J$ could be a very small subset of $\{z^{(m)}\}_{m=1}^M$ to ensure sufficient accuracy in solving the filtering problems \cite{EnSF_2023}.

\section{Numerical algorithm for the ensemble score filter in data assimilation of two-phase flow in porous media}
\label{sec:4}

In this section, we present a numerical algorithm that utilizes the training-free score-based diffusion model to formulate an ensemble score filter for solving the data assimilation problem and estimating the solution states of the two-phase flow model in porous media. Specifically, the state dynamics in the data assimilation problem are governed by a numerical solver for the two-phase flow PDE \eqref{eqn:main}. The score-based diffusion model is recursively implemented to generate samples that approximate the desired filtering distributions within the recursive Bayesian filter framework.

\subsection{Numerical solver for the two-phase flow model}\label{Sec:solver}

Here, we employ the finite element-based numerical solver for the two-phase flow as the forward dynamic model $f(\cdot)$ in the Eqn. \eqref{NLF:State}. Firstly, we reformulate Eqn. \eqref{eqn:main} as the following mixed formulation (or the first-order form):
\begin{subequations}\label{PDE_model}
    \begin{alignat}{2}
        \bu + {\bf K} \lambda(s) \nabla p &= 0, \quad \text{in } \Omega \times (0, T], \\
        \nabla \cdot \bu &= q, \quad \text{in } \Omega \times (0, T], \\
        \frac{\partial s}{\partial t} +
        \bu \cdot \nabla F(s) &= 0, \quad \text{in } \Omega \times (0, T].
    \end{alignat}
\end{subequations}
For the temporal discretization, let $s^n$, $p^n$, and $\bu^n$ denote 
the approximation of the saturation ($s(t^n)$), pressure ($p(t^n)$), and velocity ($\bu(t^n)$) at time $t^n$, where $t^n$ corresponds to the specific time at the time step $n$. Here,  $t^n = n \delta t$ with the uniform time step size $\delta t > 0$. 
%{\color{red}Given the initial conditions $\bu(\cdot, t=0) = \bu_0$ and $p(\cdot, t=0) = p_0$}

We solve the system by utilizing Implicit Pressure Explicit Saturation (IMPES) \cite{stone_analysis_1961, sheldon_one-dimensional_1959}, where pressure is solved implicitly, and then the saturation is solved explicitly.
Thus, for each time, with the given boundary conditions, saturation $s^n$, and the source term $q(t^{n+1})$, we solve for the velocity $\bu^{n+1}$ and pressure $p^{n+1}$, and then compute the saturation $s^{n+1}$. The time-discretized form is given as
\begin{subequations}
    \begin{alignat}{2}
        \bu^{n+1} + {\bf K} \lambda(s^n) \nabla p^{n+1} &= 0, \\
        \nabla \cdot \bu^{n+1} &= q(t^{n+1}), \\
        \frac{s^{n+1} - s^{n}}{\delta t} +
        \bu^{n+1} \cdot \nabla F(s^n) &= 0.
    \end{alignat}
    \label{eqn:time_discretized}
\end{subequations}

For the spatial discretization, we employ the finite element methods. We consider a shape regular, simplicial partition of the computational domain, $\overline{\Omega} = \cup_{K \in \mathcal{K}_h} \overline{K}$, where $K \in \mathcal{K}_h$ are quadrilaterals when $d = 2$. 
In this work, for the velocity, we employ the lowest order Raviart-Thomas finite element space ($\mathbb{RT}^0$), where the solution must be a vector field whose divergence is square integrable, and provides the local mass conservation. 
The spatial approximation $\bu_h$ of the velocity function $\bu(\cdot,t)$ is approximated in the following velocity finite element space:
\begin{equation}
    \mathbb{V}(K) := \left\{ \psi \in [L^2(\Omega)]^d \; \vert \; \text{div} \; \psi\vert_K \in L^2(K), \forall K \in \mathcal{K}_h \right\}.    
\end{equation}
For the pressure and saturation, we use discontinuous piecewise constant function space
\begin{equation}
    \mathbb{W}(K) := \left\{ \psi \in L^2(\Omega) \; \vert \; \psi\vert_K \in \mathbb{Q}_0(K), \forall K \in \mathcal{K}_h \right\},
\end{equation}
where $\mathbb{Q}_0(K)$ is the space of piecewise constants.

Let $\bu_h$, $p_h$, and $s_h$ be the spatial approximation of the velocity $\bu(\cdot, t)$, pressure $p(\cdot,t)$, and saturation $s(\cdot,t)$, respectively. Since we solve the velocity and the pressure before solving for saturation, we begin by discussing the bilinear for the velocity and the pressure system.  Given the initial conditions $s(\cdot, t=0) = s^0$ and the previous time step solutions $(\bu_h^n, p_h^n, s_h^n)$, we find the solution pair $(\bu_h^{n+1}, p_h^{n+1}) \in \mathbb{V}(K) \times \mathbb{W}(K)$ for $1\leq n \leq N$ such that 
\begin{subequations}
    \begin{alignat}{2}
            A(\bu_h^{n+1}, \bv_h, s_h^n) + B(\bv_h, p_h^{n+1}) &= 0, &&\quad \forall \bv_h \in \mathbb{V}(K), \\
            B(\bu_h^{n+1}, w_h) &=  -{\bf F}_v({q(t^{n+1})}, w_h), &&\quad \forall w_h \in \mathbb{W}(K),
        \end{alignat}    
\end{subequations} 
where
\begin{subequations}
\begin{alignat}{2}
A(\bu_h, \bv_h, s_h) := &\sumin\left(({\bf K} \lambda(s_h))^{-1} \bu_h, \bv_h\right)_K, \\
B({\bf v_h}, w_h) := &-\sumin(w_h, \nabla \cdot \bv_h)_K, \\
%{\color{red}{\bf F}_v(w_h)} = &\sumin(q(\cdot, t), w_h)_K.
{\bf F}_v(w_h) := &\sumin(q(t), w_h)_K,
\end{alignat}    
\end{subequations}
and $(\cdot, \cdot)_K$ denotes the $L^2(K)$-inner product on $K$.

Next, we solve for the saturation $s^{n+1}$ for given $s_h^n , \bu_h^{n+1}$, $p_h^{n+1}$, and $q_h^{n+1}$ satisfying:
\begin{equation}
Q(s_h^{n+1}, \sigma_h) = F_s(q_h^{n+1}, s_h^n, \sigma_h) - \delta t R(s_h^n, \bu_h^{n+1}, \sigma_h) \quad \forall \sigma_h \in \mathbb{W}(K),
\end{equation}
where
\begin{subequations}
\begin{alignat}{1}
Q(s_h, \sigma_h) &:= \sumin \left(s_h, \sigma_h \right)_K, \\
R(s_h, \bv_h, \sigma_h) &:= -\sumin \left( F(s_h)\bv_h, \nabla \sigma_h \right)_K + \sumin \left< F(s_h) \bn_e \cdot \bu_h, \sigma_h \right>^{up}_{\partial K}, \label{eqn:bilinear_R} \\
F_s(q_h, s_h, \sigma_h) &:= \sumin \left(s_h, \sigma_h \right) + \delta t \sumin \left( F(s_h) q_h, \sigma_h \right).
    \end{alignat}
\end{subequations}
Here, $\bn_e$ is the unit normal vector of the edge. For each quadrilateral $K \in \mathcal{K}$, the boundary integral term in Eqn. \eqref{eqn:bilinear_R} is defined in the following upwinding approach:
\begin{equation}
\left<F(s_h) \bn_e \cdot \bu_h, \sigma_h \right>^{up}_e := \left(F(s_h^+) \bn_e \cdot \bu_h^+, \sigma_h \right)_{\partial K^+} + \left(F(s_h^-) \bn_e \cdot \bu_h^-, \sigma_h\right)_{\partial K^-}, 
\end{equation}
where $\partial K^- := \{\bx \in e, \bu_h \cdot \bn_e < 0 \}$ denotes the inflow  and $\partial K^+ := \{\bx \in e, \bu_h \cdot \bn_e > 0 \}$ is the outflow. The quantities $s_h^+$ and $\bu_h^+$ denotes the values on the present element, whereas $s_h^-$ and $\bu_h^-$ denotes the quantities on the neighbor element.

\subsection{The ensemble score filter for state estimation in the two-phase flow model}

Next, we shall introduce how to use the ensemble score filter (EnSF) methodology to solve the following data assimilation problem:
\begin{equation}\label{DA:PDE}
\begin{aligned}
X_{t_{n+1}} =& \ \bar{f}(X_{t_n}, \omega_{t_n}), \qquad &\text{(State)}\\
Y_{t_{n+1}} =& \ h(X_{t_{n+1}}) + \epsilon_{t_{n+1}}, \qquad &\text{(Observation)}
\end{aligned}
\end{equation}  
where the target state $X_{t_{n}} := (s_h^{n}, \bu_h^{n}, p_h^{n})$ consists numerical solutions of the two-phase flow model computed by the state dynamics $\bar{f}(\cdot)$. 
Here,   $\bar{f}(\cdot)$ is the finite element solver for Eq. \eqref{PDE_model}, which approximates $f(\cdot) \approx \bar{f}(\cdot)$ as discussed in the previous Section \ref{Sec:solver}. The random variable $\omega_{t_n}$ in Eq. \eqref{DA:PDE} represents the uncertainty that appears in the model, which will affect the permeability and the initial conditions. 
Here, the function $h(\cdot)$ only provides partial indirect observations of $X_{t_{n+1}}$, and the observational data are perturbed by a Gaussian-type noise $\epsilon_{t_{n+1}} \sim N(0, R)$. In this setup, we recall and remark that we assume some data (e.g, ${\bf{K}}$, initial condition, solutions) are only partially observed and contain inaccuracies (noise); it serves merely as a baseline and guideline for the underlying physics in our algorithm.

The methodology of EnSF is to define score functions corresponding to the filtering distributions and approximate these scores by using Monte Carlo methods with an ensemble of state samples. Specifically, we assume that $\mathcal{S}_{n|n}$ is the approximated score function corresponding to the filtering distribution $P(X_{t_{n}}|Y_{1:t_n})$. In this way, for a temporal discretization (we assume uniform time step size 
$\Delta \tau := \tau_l - \tau_{l-1}$)
$$
0 = \tau_0 < \tau_1 < \cdots < \tau_l < \cdots < \tau_L = 1
$$
over the pseudo time interval $[0, 1]$, we can generate an ensemble of state samples, denoted by $\{z_{n|n}^{(m)}\}_{m=1}^M$, for $P(X_{t_{n}}|Y_{1:t_n})$ by solving the reverse-time SDE with the following Euler-Maruyama scheme:
\begin{equation}\label{reverse:SDE}
\begin{aligned}
\bar{Z}_{\tau_l} = \left[ b(\tau_{l+1}) \bar{Z}_{\tau_{l+1}} -\sigma^2(\tau_{l+1}) \mathcal{S}_{n|n}( \bar{Z}_{\tau_{l+1}},\tau_{l+1}) \right] \Delta \tau + \sigma(\tau_{l+1})\Delta W_{\tau}, \quad \bar{Z}_{T} \sim N(0, I_d).
\end{aligned}
\end{equation}
Then, we utilize these state samples as the initial condition at time instant $t_{n}$ and run the PDE solver $\bar{f}$ to get
\begin{equation}\label{PDE:sample}
z_{n+1|n}^{(m)} = \bar{f}(z_{n|n}^{(m)}, \omega_{t_n}^{(m)}),
\end{equation}
for sample $m = 1, 2, \cdots, M$.
The corresponding score for the prior filtering distribution $P(X_{t_{n+1}}|Y_{1:t_{n}})$ is approximated by the Monte Carlo scheme
\begin{equation}\label{prior:MC}
\bar{\mathcal{S}}_{n+1|n}(z, \tau) :=  \sum_{j=1}^{J} - \frac{z - \alpha_\tau z_{n+1|n}^{(j)}}{\beta^2_\tau} \bar{w}_\tau({z},  z_{n+1|n}^{(j)})), \quad \tau \in [0, 1]
\end{equation}
with
\begin{equation*}\label{prior:weight_app}
\bar{w}_\tau({z},   z_{n+1|n}^{(j)}) := \frac{Q_\tau(z |  z_{n+1|n}^{(j)}) }{\sum_{j=1}^{J} Q_\tau(z| z_{n+1|n}^{(j)})},
\end{equation*}
where the PDF $Q_\tau$ for the solution of the forward SDE Eq.\eqref{forward:SDE} is given explicitly by Eq. \eqref{eq:gauss}, and we name $\bar{\mathcal{S}}_{n+1|n}$ the prior score for the prior filtering distribution.

\vspace{0.5em}

To generate state samples that follow the posterior filtering distribution, we need to derive a posterior score model. Next, we describe how to update the score function and generate the posterior ensemble in EnSF by incorporating the new observational data $Y_{t_{n+1}}$. Specifically, we update the prior score $S_{n+1|1:n}$ with the new observation $Y_{t_{n+1}}$ according to the Bayes's rule from Eq.~\eqref{Update} to approximate the posterior score $S_{n+1|1:n+1}$:
\begin{equation}\label{eq:EnSF-posterior-score}
    \bar{S}_{n+1|1:n+1}(z, \tau) = \bar{S}_{n+1|1:n}(z, \tau) +  h(\tau) \nabla_{z} \log P(Y_{t_{n+1}} | X_{t_{n+1}})(z),
\end{equation}
where $\bar{S}_{n+1|n}$ is the Monte Carlo estimator of the prior score in Eq.~\eqref{eq:score11}, $P(Y_{t_{n+1}} | X_{t_{n+1}})(\cdot)$ is the likelihood function in Eq.~\eqref{Update}, and $h$ is a continuous time-damping function to control the diffusion of the information from $Y_{t_{n+1}}$ in the diffusion domain. In the current EnSF method, the function $h(t)$ is a monotonically decreasing in $[0,1]$ (e.g., $h(\tau) = 1-\tau$ in the numerical experiments in this paper) satisfying $h(1) = 0$ and $h(0) = 1$, which indicates that the information in the likelihood is gradually incorporated into the score function while solving the reverse SDE.

\subsection{Summary of the  data assimilation framework for the two-phase flow model.}

In this subsection, we summarize the overall data assimilation (DA) framework for the two-phase flow model.  

%At each time step $t = n$, given the state vector $X_{t_n}$, we use the Finite Element Method (FEM) solver to obtain the prior prediction $\hat{X}_{t_{n+1}}$ for the DA method. 

\begin{itemize}
    % \item \textit{Construct reference data:} \\
    % For all, $t_n \in (0, T]$, where $n=1, \cdots, N$, we obtain the solution set
    % $$
    % \hat{\mathbb{X}} = \{ \hat{X}_{t_0}, \cdots, \hat{X}_{t_N} \}
    % $$
    % by solving Eq.\eqref{eqn:time_discretized} using the Finite Element Method (FEM), where the permeability field \textbf{K} may have spatial noise but remains constant over time. The sequence of data is generated using 
    % $$
    % \hat{X}_{t_{n}} = \bar{f}(\hat{X}_{t_{n-1}}, \textbf{K}),
    % $$
    % where $\hat{X}_{t_0}$ comes from the initial condition. 
    % %where $\bar{f}(\cdot)$ is the FEM solver. %where ${\color{red}\omega_{t_{n}}}$ is the input to perturb the data such as the permeability, i.e  $K({\color{red}\omega_{t_{n}}})$. 
    % %{\color{red} In this step, we do not  perturb the permeability, and we consider the solution set $\tilde{\mathbb{X}}$ as the real observational data.}
    \item For each time step $t_n$,
    \begin{enumerate}
        \item %\textit{FEM step:} $X_{t_{n+1}} =  \bar{f} (\tilde{X}_{t_{n}}, \omega_{t_n})$
        When $n=0$, the initial condition of PDE is given as $\tilde{X}_{t_{0}}$, and the initial score $\bar{\mathcal{S}}_{n|n}$ is given corresponding to $\tilde{X}_{t_{0}}$. %, i.e.,  $\tilde{X}_{t_{0}} = \hat{X}_{t_0} + \eta$.
        % \begin{itemize}
        %     %\item if $n=0$, the initial condition is perturbed with a noise $\eta$, i.e.,  $\tilde{X}_{t_{0}} = \hat{X}_{t_0} + \eta$. The uncertainty in the solution $X_{t_{n+1}}$ arises from both the initial condition and random variable $\omega_{t_n}$, which perturbs the permeability. 
        %     %\item if $n \geq 1$, $\tilde{X}_{t_{n}}$ is computed by the following DA algorithm 
        % \end{itemize}
        \item For $n \geq 0$, implement the data assimilation framework to obtain the estimated PDE solution:
        $$
        {\tilde{X}_{t_{n+1}}} = \text{DA}(Y_{t_{n+1}}, {\tilde{X}_{t_{n}}} ),
        $$
        where $Y_{t_{n+1}}$ is the partial noisy observation of the true $X_{t_{n+1}}$.
        
        %{\color{red}where the details of $\text{DA}$ algorithm is given in Algorithm \ref{alg:posterior_sampling}}.
        
        \textit{Elaboration of the DA$(Y_{t_{n+1}}, {\tilde{X}_{t_{n}}} )$ operator}
    \begin{enumerate}
    \item We set the number of ensembles $M$, the number of pseudo timesteps $L$ required in the reverse SDE in the ensemble score filter.
    \item For $\bar{\mathcal{S}}_{n|n}$ from the previous time step corresponding to $p(\tilde{X}_{t_n}|Y_{t_n})$, we use this score function to generate the state samples $\{z_{n|n}^{(m)}\}_{m=1}^{M}$ by using equation \eqref{reverse:SDE}. The samples $\{z_{n|n}^{(m)}\}_{m=1}^{M}$ correspond to $\tilde{X}_{t_{n}}$. 
    \item For each of the state sample $\{z_{n|n}^{(m)}\}_{m=1}^M$, we use the FEM solver $\bar{f}$ (forward solver) to generate the prior prediction $\{z_{n+1|n}^{(m)}\}_{m=1}^{M}$ and we use $\eqref{prior:MC}$ to estimate $\bar{\mathcal{S}}_{n+1|n}$
    \item Update with scheme \eqref{eq:EnSF-posterior-score} to obtain $\bar{\mathcal{S}}_{n+1|n+1}$.
    \item Generate samples $\{z_{n+1|n+1}^{(m)}\}_{m=1}^M$ from $\bar{\mathcal{S}}_{n+1|n+1}$ to obtain the DA estimated solution ${\tilde{X}_{t_{n+1}}}$.
    \item[] \textbf{Remark.} In our numerical experiments, we set $M=300$ for each time steps, and $L=1000$.
\end{enumerate}

    \end{enumerate}
    %\item \textbf{FEM Step ($n=1$)}. For the first (initial) time step ,  with the given initial data 
    %$\tilde{X}_{t_0}, $
    %we obtain $${X}_{t_{1}} = \text{FEM}(X_{t_0} , \bar{K}(\omega)).$$
\end{itemize}

We emphasize that the permeability used to generate the reference solution for observational data differs from that used in the FEM solver for prediction within the data assimilation framework. This discrepancy reflects a common challenge in real-world applications: incomplete observations of key model parameters and the presence of intrinsic model errors and uncertainties in numerical simulations.

% Throughout alg.1, we keep M ensembles and do forward simulation for each of it, therefore the $\hat{X}_{t_{n}} = \{z_{n|n}^{(m)}\}_{m=1}^M$ for each $m = 1,2,...M$. All $z_{n|n}^{(m)}$ is generated from the reverse-time SDE in \eqref{reverse:SDE} from $\bar{Z}_{T} \sim N(0, I_d)$ . When n = N s.t. we reach the terminal time, we average $\{z_{n|n}^{(m)}\}_{m=1}^M$ for all $M$ ensembles such that we have a single output.

\section{Numerical experiments}
\label{sec:5}
In this section, we present three numerical examples to evaluate the performance of our algorithm. A key motivation for this research stems from the fact that, in real-world applications, the numerical PDE model used in simulations often differs from the ``perfect'' physical model due to limitations in our knowledge. To reflect this reality in our experiments, we initialize the PDE solvers with different initial conditions and assume \textit{limited knowledge} of the true permeability when solving the two-phase flow model using finite element solvers. 

In the numerical experiments section, all the experiments were implemented on a 16-core Intel(R) Core(TM) i7-10700 CPU @ 2.90GHz and NVIDIA RTX 2060. 

\subsection{Experimental Setups}
\subsubsection{FEM solver}
Here, we describe the setup of the FEM solver for the following numerical experiments. The computational domain is given as $\Omega = [0,1]^2$ for $t \in (0,T]$. The total  mobility is defined as 
\begin{equation}
    \lambda (s) = \frac{1.0}{\mu}s^2 + (1-s)^2, 
\end{equation}
where the viscosity $\mu = 0.2$~\cite{2024:africa.arndt.ea:deal}. In addition, the fractional flow is given by 
\begin{equation}
    F(s) := \frac{s^2}{s^2 + \mu (1-s)^2}.
\end{equation}
We assume that there is no source term, $q := 0$, and the porous medium is isotropic,
\begin{equation}
    \mathbf{K}(\bx) = k(\bx)\mathbf{I}, 
    \label{eqn:isotropic_permeability}
\end{equation}
where $\mathbf{I}$ is the identity matrix. 
%{\color{red}For example 4, the permeability is shown in Figure X.}

To supplement the governing system, the following  initial and boundary conditions are imposed: 
\begin{subequations}
    \begin{alignat}{2}
        p(\bx, t) &= 1 - x_1 \quad && \text{on} \; \partial \Omega, \\
        s(\bx, t) &= 1 \quad && \text{on} \; \Gamma_{in} \in {x_1 = 0}, \\
        s(\bx, t) &= 0 \quad && \text{on} \; \Gamma_{in} \textbackslash {x_1 = 0}, \\
        s(\bx, 0) &= 0 \quad && \text{on} \; \Omega,       
    \end{alignat}
\end{subequations}
where $\Gamma_{in}(t) = \{\bx \in \partial \Omega :  \bu(\bx, t) \cdot \bn< 0\}$.

\subsubsection{Setup for DA algorithm (EnSF)}
{Due to the limitation of getting the real observational data, we use synthetic data as our reference data to construct the observation data. To generate the reference data, for all $t_n \in (0, T]$, where $n=1, \cdots, N$, we obtain the solution set
$$ \hat{\mathbb{X}} = \{ \hat{X}_{t_0}, \cdots, \hat{X}_{t_N} \} $$ by solving Eq.\eqref{eqn:time_discretized} using the Finite Element Method (FEM), where the permeability field \textbf{K} may have spatial noise but remains constant over time. The sequence of data is generated using 
$$ \hat{X}_{t_{n}} = \text{FEM}(\hat{X}_{t_{n-1}}, \textbf{K}), $$
where $\hat{X}_{t_0}$ comes from the initial condition and $\text{FEM}(\cdot, \cdot)$ is the finite element solver for Eqn.\eqref{eqn:time_discretized} with the boundary conditions.} \par
In the numerical examples that follow, we incorporate mild uncertainty in the initial condition. Specifically, the initial velocity and pressure are set to zero, while the saturation is perturbed by small Gaussian noise $N(0, 0.01 I_d)$. The primary source of uncertainty in our study lies in the permeability term $\bf{K}$. Accurately estimating permeability is notoriously difficult in practical geophysical modeling, and its parametrization remains a significant challenge.

To solve the PDE model numerically, we discretized the domain with a $64 \times 64$ mesh in space. This makes the EnSF solving a {$16,512$}-dimensional data assimilation problem -- with $4,096$ dimensions for the saturation $s_h$, {$8,320$} dimensions for the velocity $\bu_h$ and $4,096$ dimensions for the pressure $p_h$. 

%{\color{red} Throughout the numerical examples, we focus on the saturation values. The observation data and the saturation solution from $\bar{f}(\cdot)$ at a given time $t_n$ are denoted as $\hat{X}_{t_n} = \hat{s}_h(t_n)$, and ${X}_{t_n} = {s}_h(t_n)$, respectively. We note that ${X}_{t_n} = {s}_h(t_n)$ is the numerical solution directly given from $\bar{f}(\cdot)$, without any data assimilation. Then, $\tilde{X}_{t_n} = \tilde{s}_h(t_n) $ is the final solution that we obtain with the proposed EnSF data assimilation algorithm. }

We want to point out that due to the curse of dimensionality issue for approximating high-dimensional distributions, a {$16,512$}-dimensional data assimilation problem is extremely challenging since one needs to approximate the filtering distribution $P(X_{t_n}|Y_{t_1:t_{n}})$ in the {$16,512$}-dimensional state space. We also want to mention here that by partial observation, for example, $50\%$ observation in saturation means that we only observe randomly picked $50\%$ area for saturation information. For the unobserved area, they do not provide observational information.

\subsection{Example 1. A case with uncertainty in permeability matrix}

In real-world scenarios, the true permeability of the subsurface may differ substantially from the values assumed in the model. This might result in model errors that could evolve unpredictably over time. In this example, we add uncertainties in the permeability matrix to represent our incomplete knowledge in formulating the two-phase flow model. \par
\begin{figure}[h!]
    \centering
    \begin{subfigure}[t]{0.3\textwidth}
        \includegraphics[width=\textwidth]{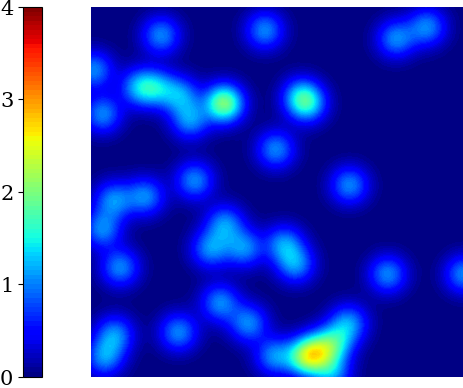}
    \caption{$\hat{k}(\bx)$}
    \label{Ex1:k_hat}
    \end{subfigure} 
    \hspace{1cm}
    \begin{subfigure}[t]{0.3\textwidth}
        \includegraphics[width=\textwidth]{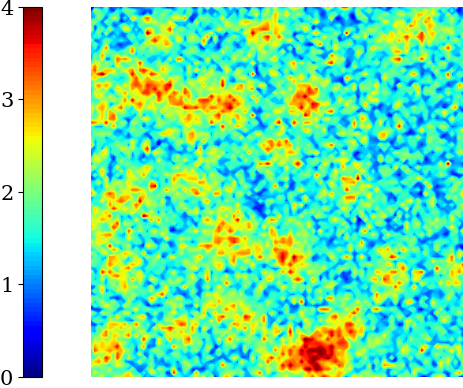}
    \caption{$k(\bx)$}
    \label{Ex1:k}
    \end{subfigure}
    \caption{Permeability function}
    \label{Ex1:permeability}
\end{figure}

Here, we introduce a time-invariant uncertainty that affects 80\% of the spatial domain. To construct the reference data $\hat{X}$ for all t, the permeability matrix $\textbf{K}(\bx)$ (as given in Eq.\eqref{eqn:isotropic_permeability}) is computed using the following $k(\bx)$:
\begin{equation}\label{Ex1:K}
k(\bx) = \min\left\{\hat{k}(\bx) + \eta_1(\bx) + \eta_2(\bx) + \eta_3(\bx), 4\right\},
\end{equation}
where $\{\eta_i(\bx)\}_{i=1}^3$ represent the uncertainties defined as follows:
\begin{itemize}
\item $\eta_1$ occurs in $56\%$ of the region, following $\eta_1(\bx) \sim N(1.5, 0.75)$, 
\item $\eta_2$ occurs in $16\%$ of the region, following $\eta_2(\bx) \sim N(1.0, 0.5)$,
\item $\eta_3$ occurs in $8\%$ of the region, following $\eta_3(\bx) \sim N(2.0, 1.0)$.
\end{itemize}
The base permeability function $\hat{k}(\bx)$ is defined as
\begin{equation}\label{def:K}
    \hat{k}(\bx) := \min \left\{ \max \left( \sum_{i=1}^{N} \sigma_i(\bx), 0.01 \right), 4 \right\}, \ \
    \sigma_i(\bx) := \exp\left({-400\left( {\vert \bx - \bx_i \vert} \right)^2}\right),
\end{equation}
where
$\{\bx_i\}$ are $N$ randomly chosen locations inside the domain, which models a domain with $N$ centers of higher permeability regions \cite{2024:africa.arndt.ea:deal}. Fig. \ref{Ex1:permeability} shows the base permeability function, $\hat{k}(\bx)$, alongside the permeability with added noise, $k(\bx)$. Fig. \ref{Ex1:k_hat} is the permeability obtained using $N=40$ randomly selected centers. Fig. \ref{Ex1:k} shows the permeability, $k(\bx)$, after adding the Gaussian noise to 80\% of the spatial domain, as described in Eq.\eqref{Ex1:K}.\par
After $\{\eta_i(\bx)\}_{i=1}^3$ are generated, the randomness in permeability is fixed. Moreover, we assume no prior knowledge of these uncertainties in the numerical experiments to examine the robustness of the EnSF method. Moreover, we provide indirect and noisy observations to our data assimilation algorithm. We define the observation as
\begin{equation}\label{NL}
    Y(\hat{X}_{t_{n+1}}) := Mask(arctan(\hat{X}_{t_{n+1}}) + \epsilon),
\end{equation}
where $Mask(\cdot)$ is a selection operator that simulates the situation when detectors are placed randomly to observe some data including saturation values, velocity, and pressure. In this example, we only observe the saturation and do not apply data assimilation to the velocity and the pressure. 
The data is observed through the $arctan(\cdot)$ observational function, and $\epsilon \sim \mathcal{N}(0,0.07 I_{d})$ adds Gaussian noises to the observed states.

\begin{figure}[h!]
    \centering
     %\hfill
    \begin{subfigure}[t]{0.3\textwidth}
        \includegraphics[width=\textwidth]{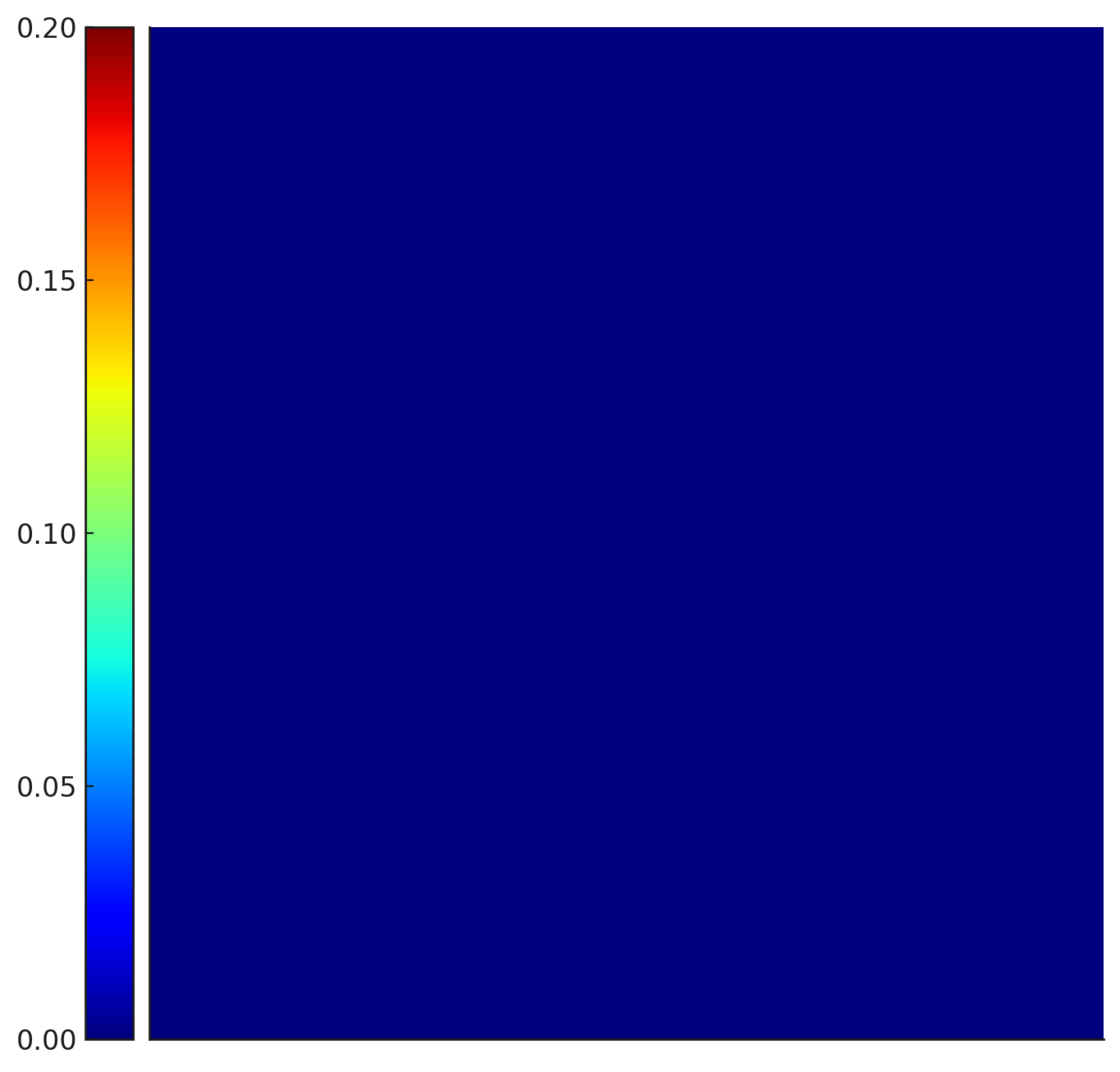}
    \caption{Initial condition for the reference saturation, $\hat{s}_h(\cdot, t=0)$}
    \label{Ex1:IC_ref_Saturation}
    \end{subfigure} 
    \begin{subfigure}[t]{0.3\textwidth}
        \includegraphics[width=\textwidth]{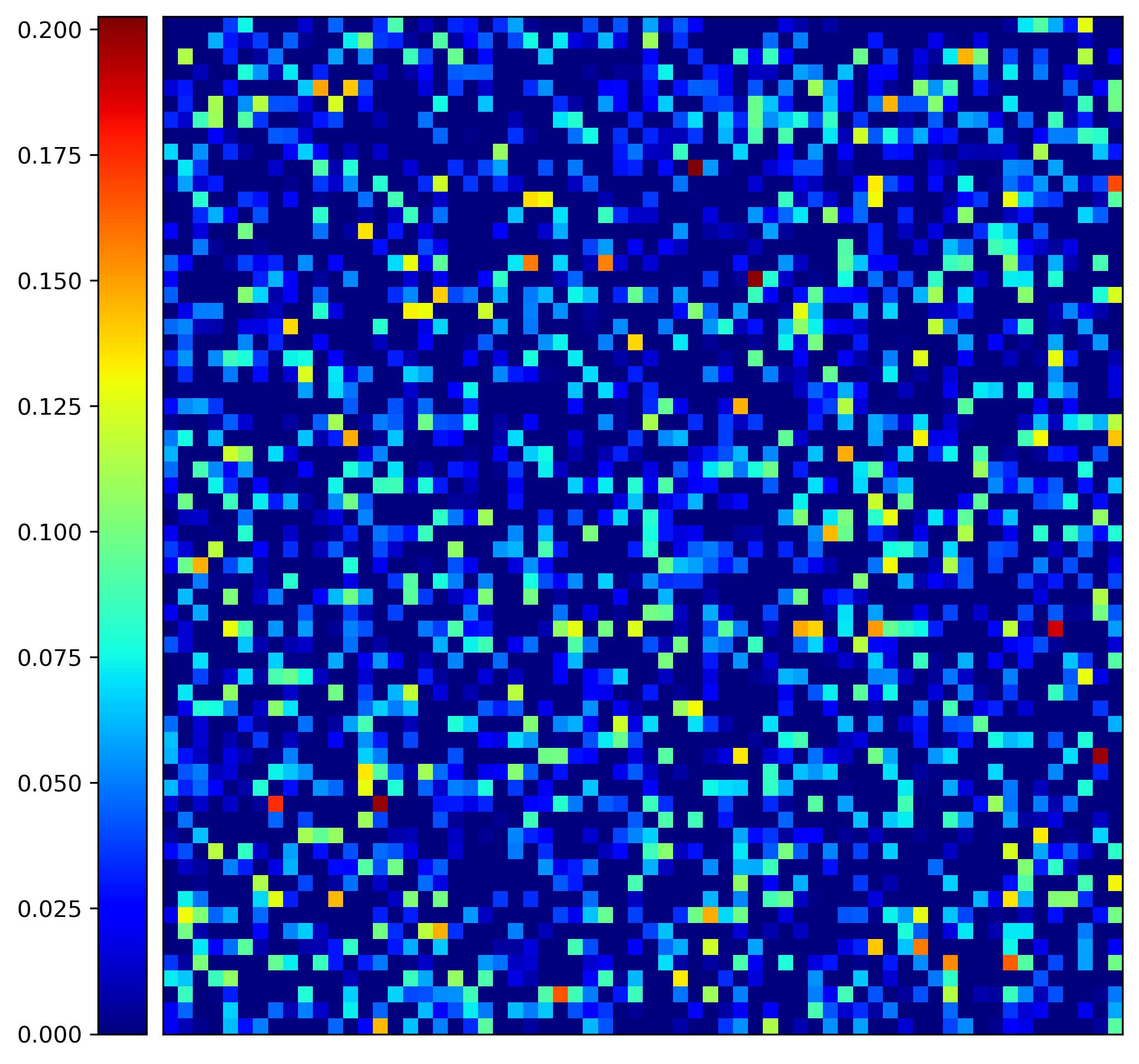}
    \caption{Initial condition for the saturation $s_h$ in the forward model $\bar{f}$}
    \label{Ex1:IC_Saturation}
    \end{subfigure} 
    \caption{Initial conditions of the saturation. a) illustrates the initial condition for the reference saturation $\hat{s}_h(t=0) = 0 $ and b) presents the perturb initial condition for solving $s_h$ in the forward model $\bar{f}$.}
    \label{Ex1:Saturation}
\end{figure}

Thus, in our DA algorithm, we use $\hat{k}(\bx)$ from Eq. \eqref{def:K} in our forward dynamic model $\bar{f}(\cdot)$ and use partial observations that are noisy and indirect. 
%\textcolor{cyan}{QUESTION: IN THE ALGORITHM IN SECTION 4.3, WE WROTE $\bar{f}$ HAS A RANDOM VARIABLE $\omega_{t_n}$. WHAT IS THIS RANDOM VARIABLE IN THIS EXAMPLE?} \textcolor{purple}{To be deleted: This $\epsilon$ comes from the observation equation in \eqref{DA:PDE} for the second equation. Basically it is also a random variable but independent of the $\omega$ in model. This is a usual description in data assimilation to model the scenario that we know things we observe is also not going to be completely accurate due to the measurement inaccuracy. A good example would be considering measuring the distance between two objects by using some facility that can only reads one decimal digit, if the meter reads 1.5 meter, it doesn't necessarily mean that two objects are 1.5m away from each other. It might be 1.54xxxx meters away from each other but starting from the second digit, things are dropped. In this case the observation is direct and $\epsilon = 0.04xxxx$.  }  
Furthermore, since our knowledge of the domain's initial condition may be limited, we address this uncertainty by perturbing the initial condition for saturation. We use a random initial condition that follows the absolute value of the standard normal distribution, i.e, 
\begin{equation}
    \tilde{s}_h^0  \sim \vert N(0, 1/300) \vert.
\end{equation}
The initial condition for the reference model and forward model is shown as Fig. \ref{Ex1:Saturation}\par
\begin{figure}[h!]
    \centering
        \begin{subfigure}[t]{0.3\textwidth}
            \centering
            \includegraphics[width=\textwidth]{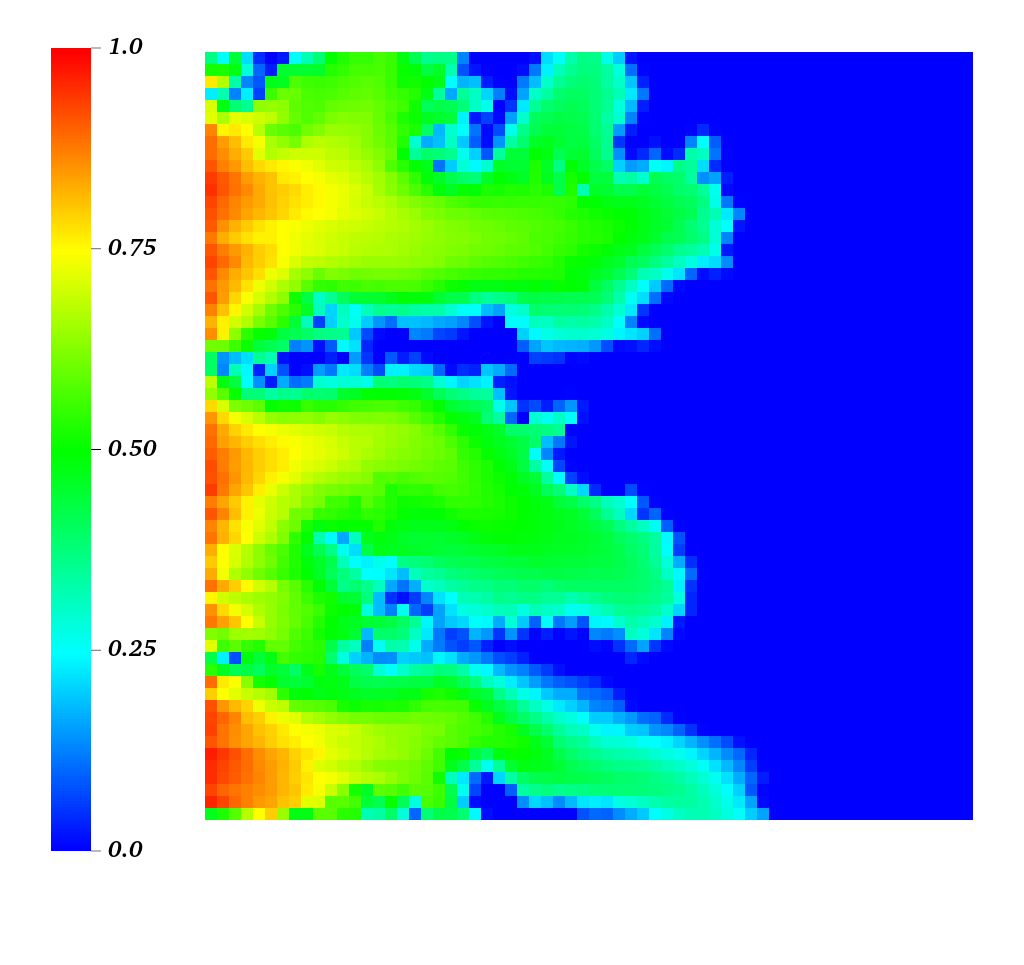}
            \caption{Reference saturation $(\hat{s}_h)$}
            \label{Ex1:Refa}
        \end{subfigure} 
        \begin{subfigure}[t]{0.3\textwidth}
            \centering
            \includegraphics[width=\textwidth]{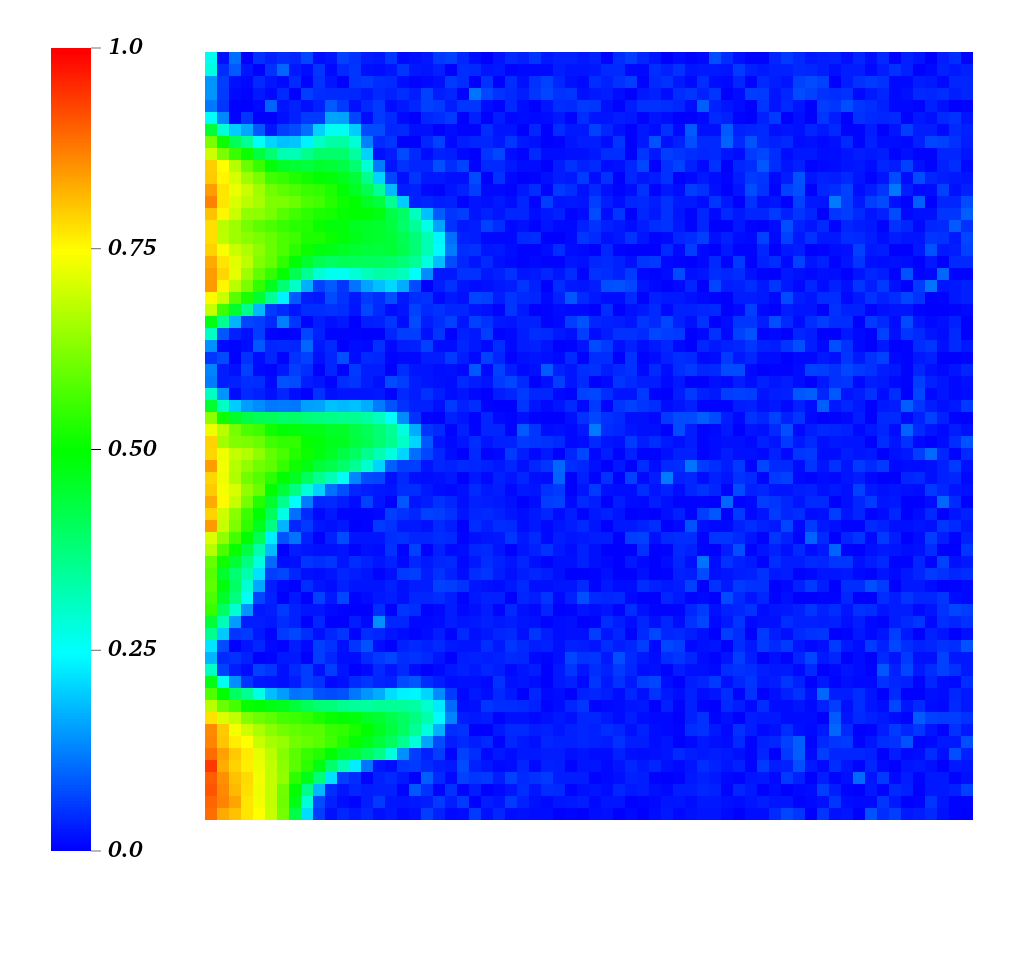}
            \caption{Estimated saturation with no observation, i.e., the saturation computed solely by the state dynamics $f(\cdot)$ using $\hat{k}(\bx)$ as the permeability}
            \label{Ex1:Refb}
        \end{subfigure} 
        \begin{subfigure}[t]{0.3\textwidth}
            \centering
            \includegraphics[width=\textwidth]{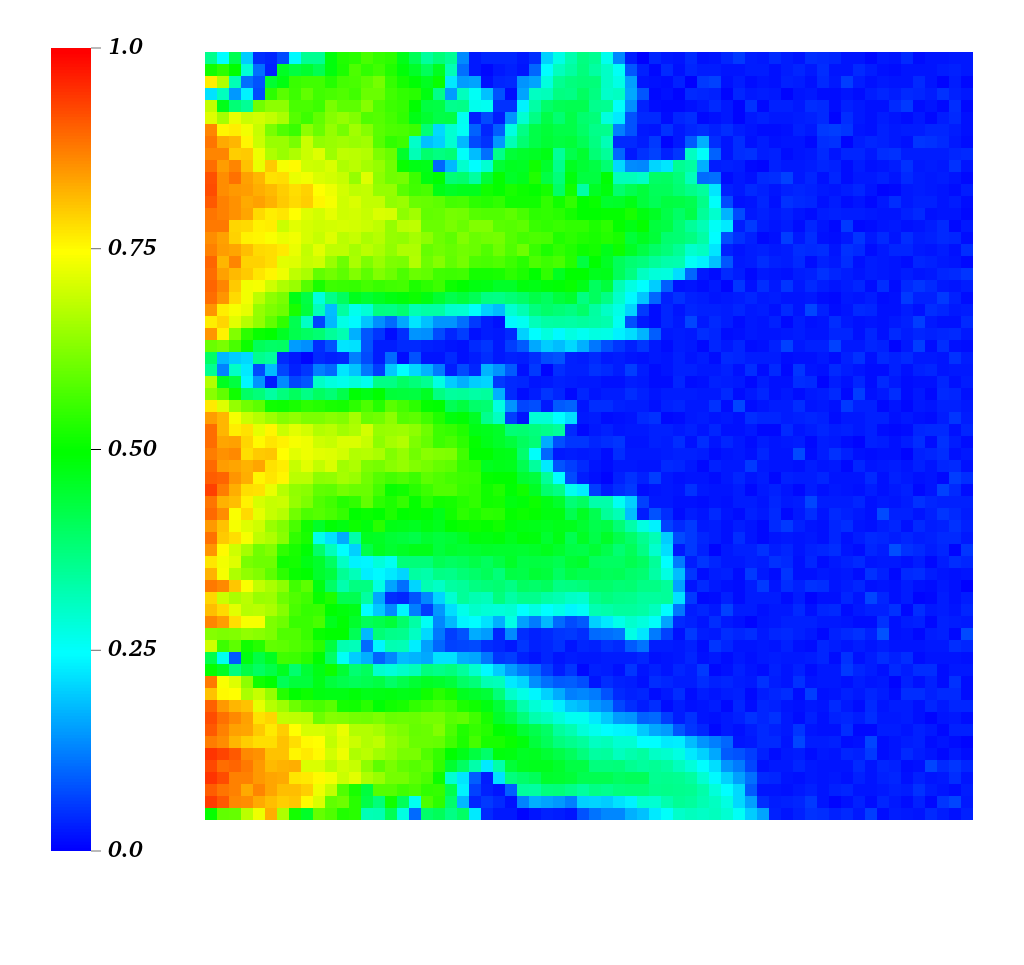}
            \caption{Estimated saturation $(\tilde{s}_h)$ with full observation of the saturation }
            \label{Ex1:Refc}
        \end{subfigure} 
    \caption{The comparison between the reference saturation, estimated saturation without observation, and the estimated saturation with EnSF at $t = 0.4$. a) the reference saturation $(\hat{s}_h)$ is simulated using $k(\bx)$ and the true initial condition, without any perturbation. 
    b)-c) Both estimated saturations are computed from noisy initial saturation-- b) one without observation, and c) the other using the EnSF with 100\% observation of the saturation data.}
    \label{Ex1:Ref}
\end{figure}

In Fig. \ref{Ex1:Ref}, we compare the results with and without our data assimilation algorithm, using a time step of $\delta t = 0.001$. Fig. \ref{Ex1:Refa} shows the reference data, saturation ($\hat{s}_h$), at $t = 0.4$. Fig. \ref{Ex1:Refb} is the final saturation obtained without the data assimilation, i.e., the saturation at $t = 0.4$ by the state dynamics $\bar{f}(\cdot)$ using $\hat{k}(\bx)$ with a random initial condition. Without the data assimilation, we observe a significant deviation from the reference saturation (Fig. \ref{Ex1:Refa}), as no data assimilation procedure is applied to calibrate the simulated solution obtained by the forward model. On the other hand, using the EnSF-based data assimilation with the full observation of the saturation, we recover the essential solution features. Fig. \ref{Ex1:Refc} shows that the result is in good agreement with the reference solution. This highlights the effectiveness of our algorithm in recovering the underlying true solution, even in the presence of uncertainty of the permeability and noisy initial conditions. \par
\begin{figure}[h!]
    \centering
        \begin{subfigure}[t]{0.3\textwidth}
            \centering
            \includegraphics[width=\textwidth]{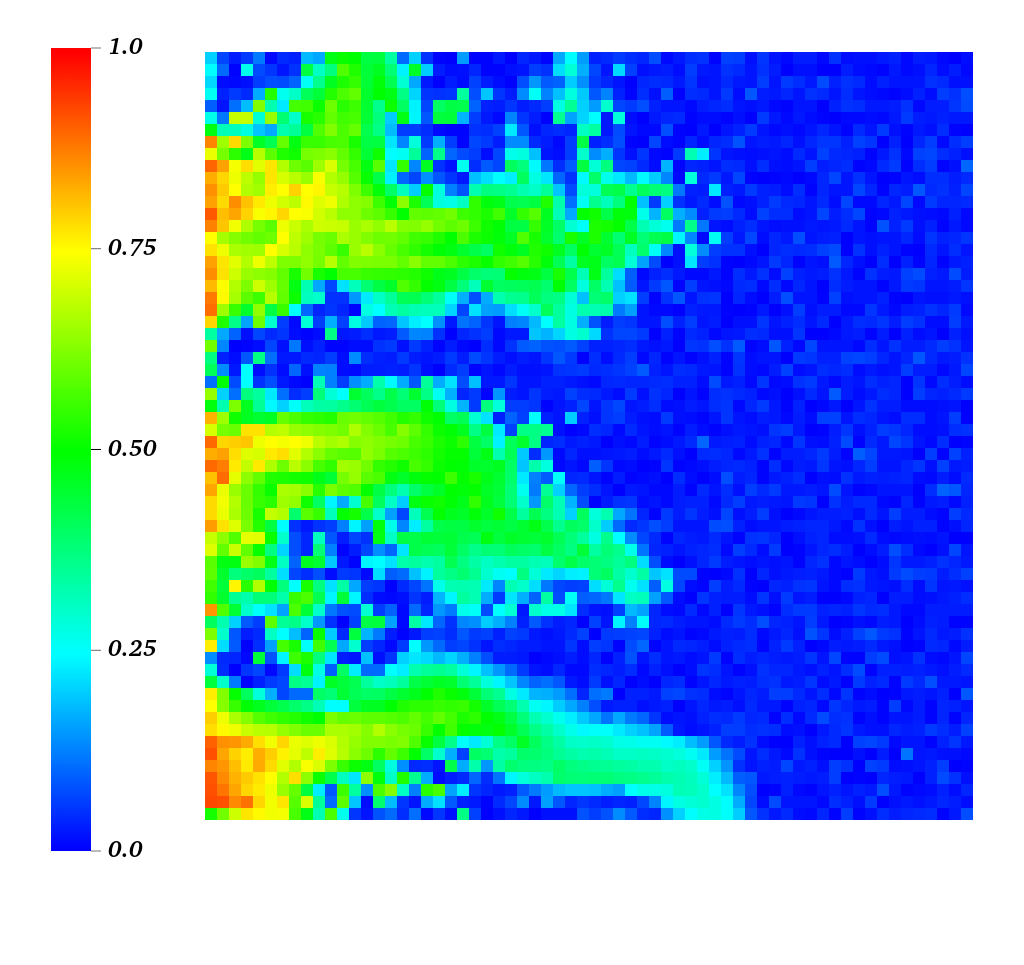}
            \caption{$25\%$ observations}
            \label{Ex1:Partial_a}
        \end{subfigure}
        \begin{subfigure}[t]{0.3\textwidth}
            \centering
            \includegraphics[width=\textwidth]{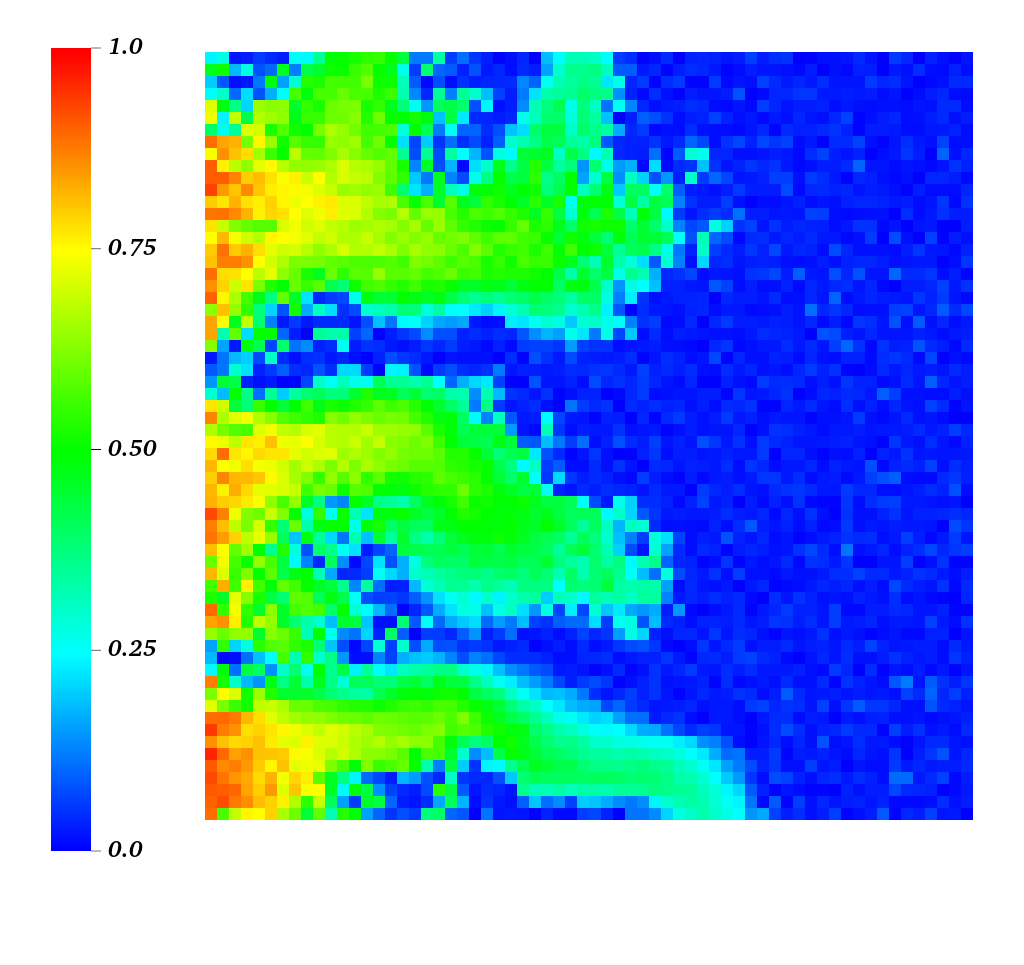}
            \caption{$50\%$ observations}
            \label{Ex1:Partial_b}
        \end{subfigure} 
        \begin{subfigure}[t]{0.3\textwidth}
            \centering
            \includegraphics[width=\textwidth]{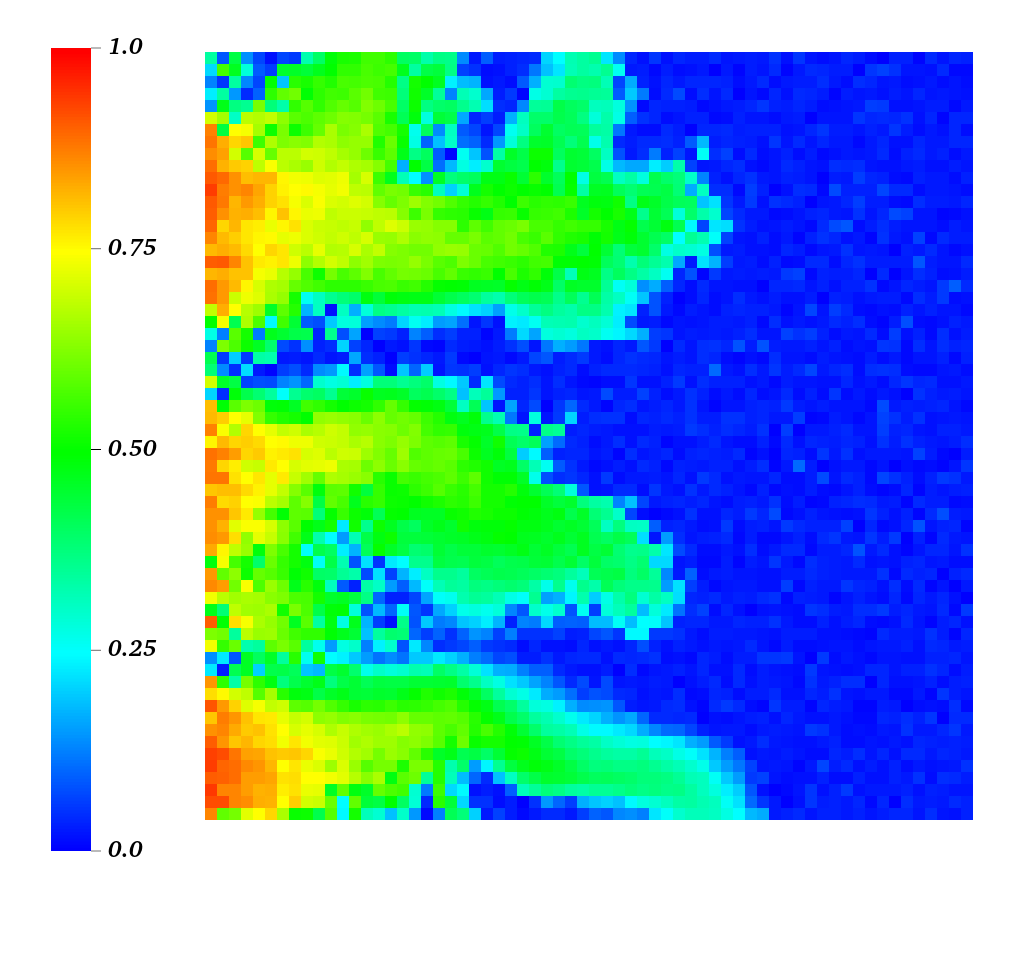}
            \caption{$75\%$ observations}
            \label{Ex1:Partial_c}
        \end{subfigure} 
    \caption{Estimated saturation ($\tilde{s}_h$) with EnSF based on partial observations of the saturation data.}
    \label{Ex1:Partial}
\end{figure}

\begin{figure}[h!]
    \centering
        \begin{subfigure}[t]{0.3\textwidth}
            \centering
            \includegraphics[width=\textwidth]{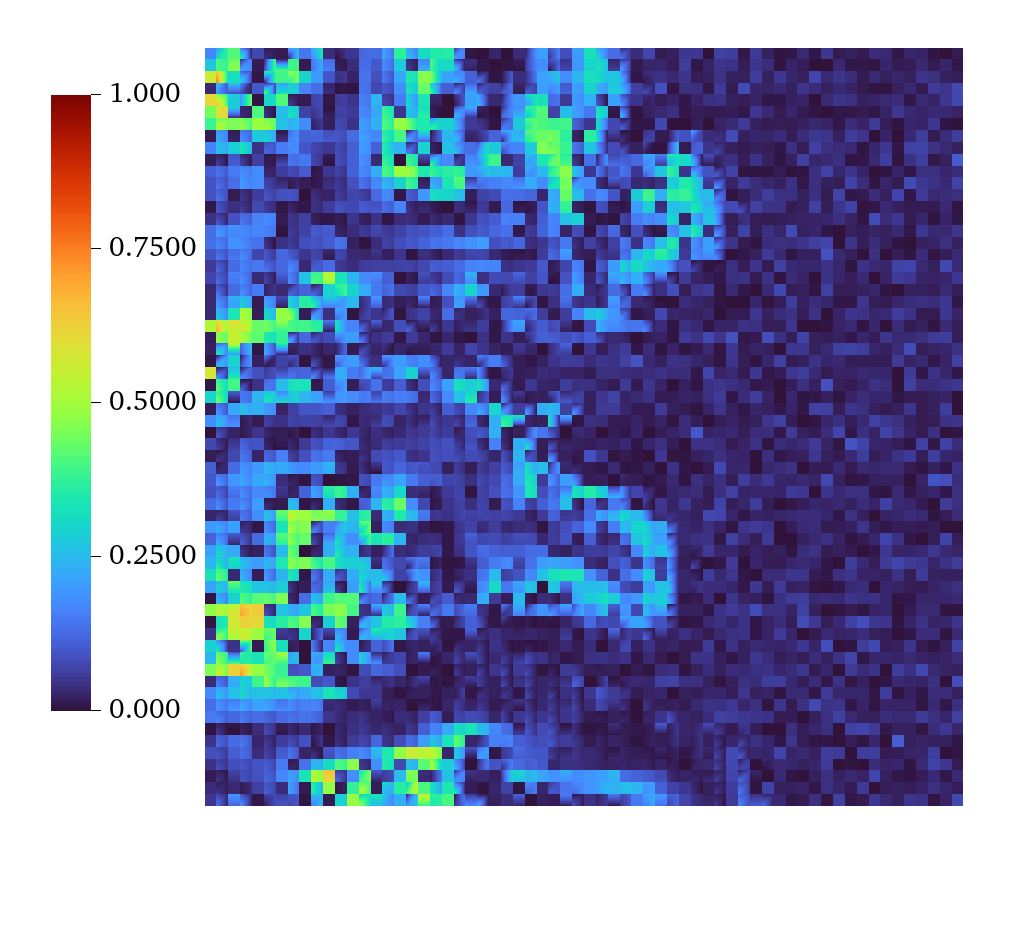}
            \caption{$25\%$ observations}
            \label{Ex1:Error_a}
        \end{subfigure}
        \begin{subfigure}[t]{0.3\textwidth}
            \centering
            \includegraphics[width=\textwidth]{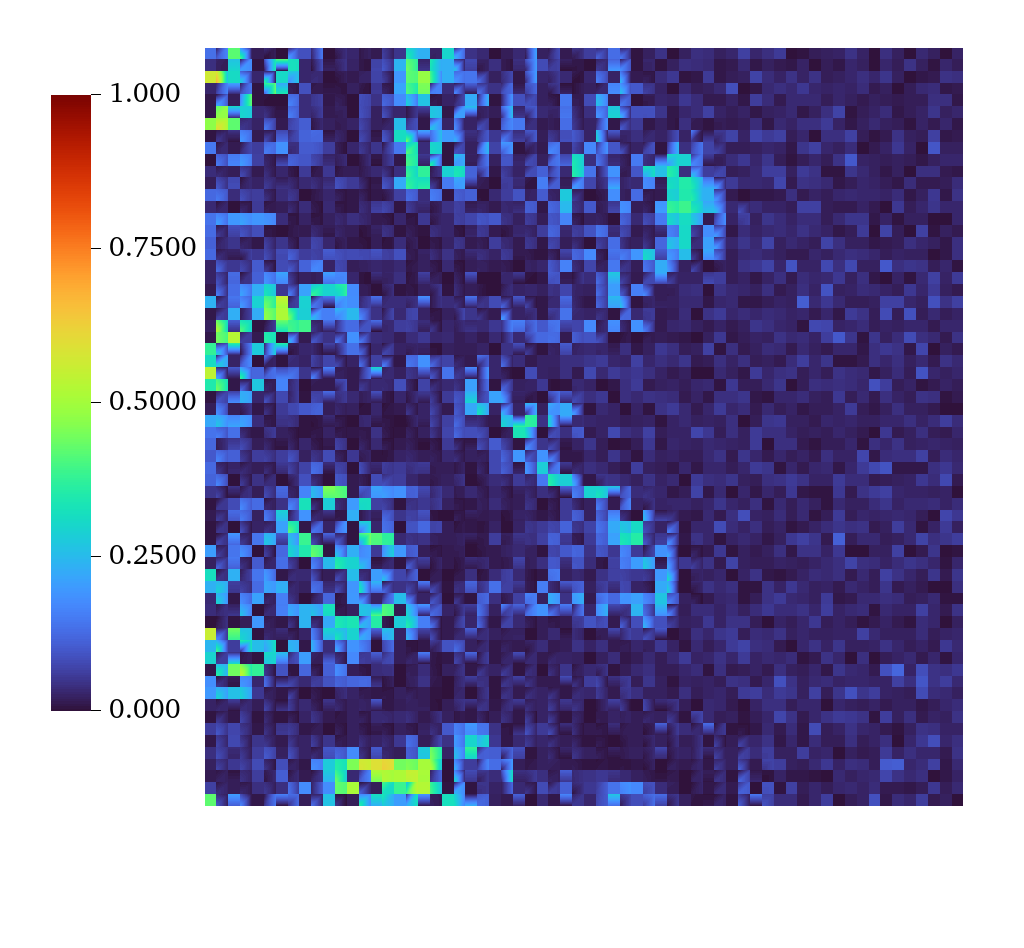}
            \caption{$50\%$ observations}
            \label{Ex1:Error_b}
        \end{subfigure} 
        \begin{subfigure}[t]{0.3\textwidth}
            \centering
            \includegraphics[width=\textwidth]{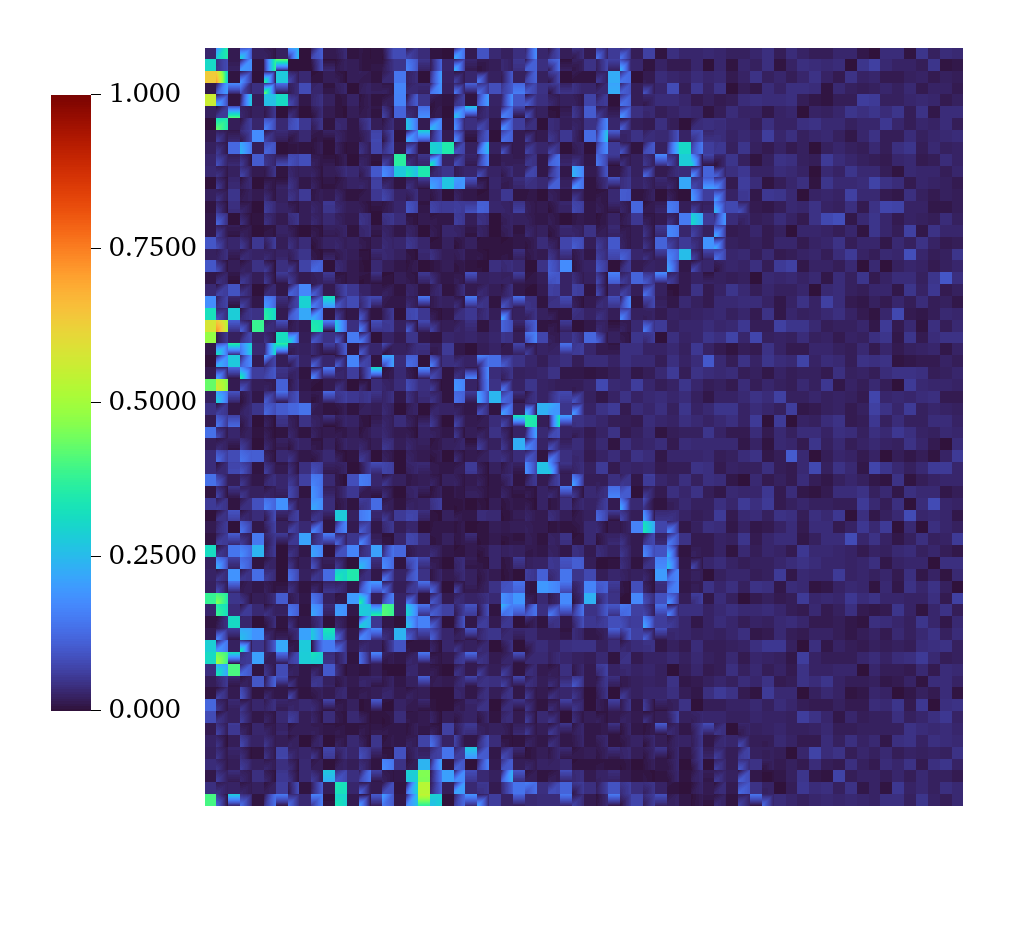}
            \caption{$75\%$ observations}
            \label{Ex1:Partial_c}
        \end{subfigure} 
    \caption{Error for partial observations of the saturation data.}
    \label{Ex1:Error}
\end{figure}

\begin{figure}[ht]
    \centering
    \begin{subfigure}[t]{0.5\textwidth}
        \includegraphics[width=\textwidth]{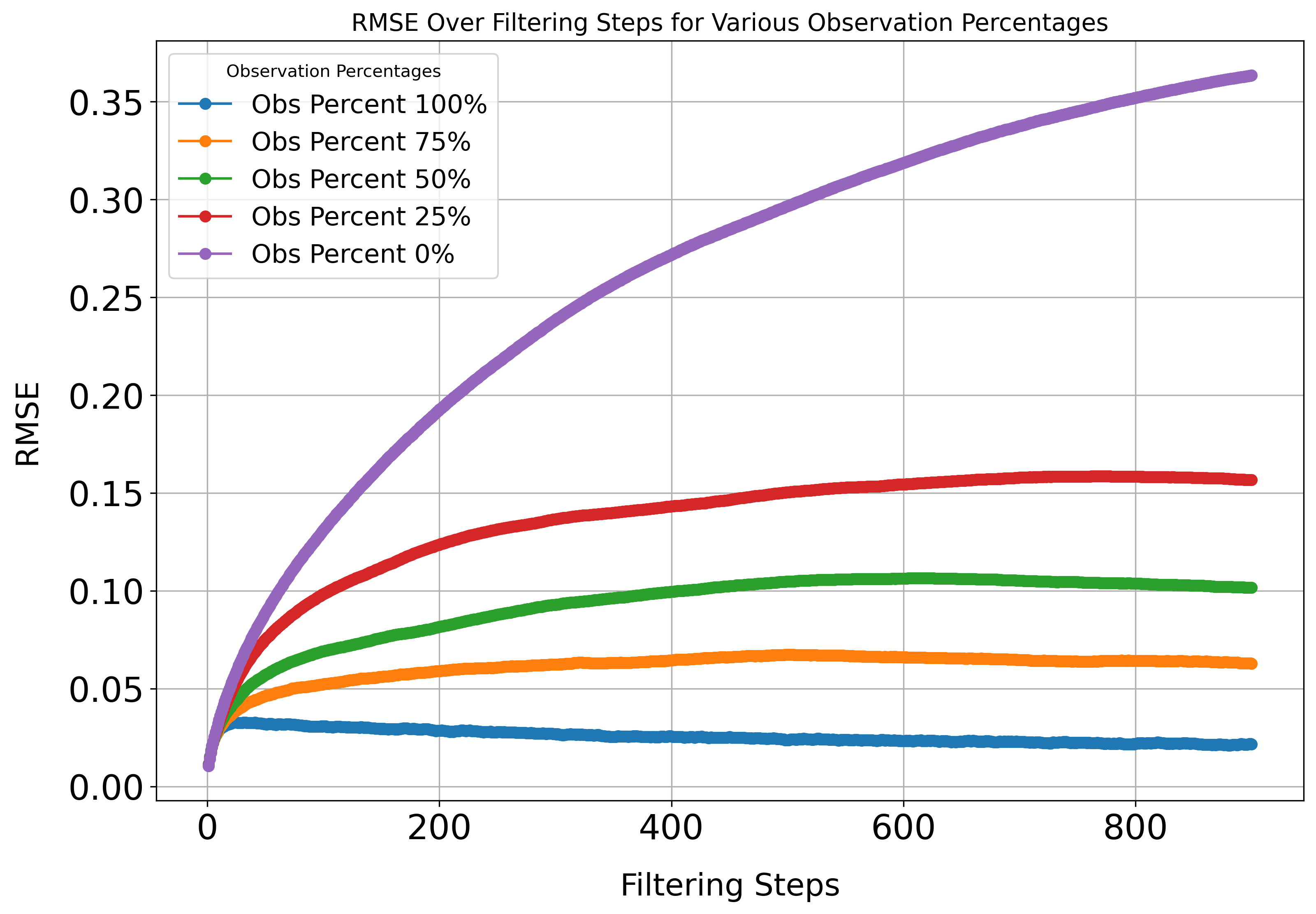}
    \end{subfigure} 
    \caption{RMSEs comparisons of 0\%, 25\%, 50\%, 75\% and 100\% observations}
    \label{Ex1:RMSE}
\end{figure}

Due to limited resources in the real world, we cannot observe the situation across the domain and typically only have access to partial observations. Here, we highlight that our EnSF algorithm is capable of recovering the underlying solution only with partial observations. Fig. \ref{Ex1:Partial} compares the performance of the algorithm across various levels of observations of the saturation. We can see in Fig. \ref{Ex1:Partial_a} that even with only 25\% of the observation of the saturation, the EnSF effectively captures the main features of the saturation, and incorporating more data leads to increasingly accurate estimates of the solution. Fig. \ref{Ex1:Error} shows the absolute difference between the reference saturation ($\hat{s}_h$) and the estimated saturation ($\tilde{s}_h$) for various levels of observations. As the observation increases, the error in the estimated saturation decreases. In Fig. \ref{Ex1:RMSE}, we show the root mean square errors (RMSEs) over the filtering steps for different observation percentages. It shows that with an increasing number of observations, our method provides better calibration of the result compared to the reference solution. 

\begin{figure}[h!]
    \centering
    \begin{subfigure}{0.3\textwidth}
        \includegraphics[width=\linewidth]{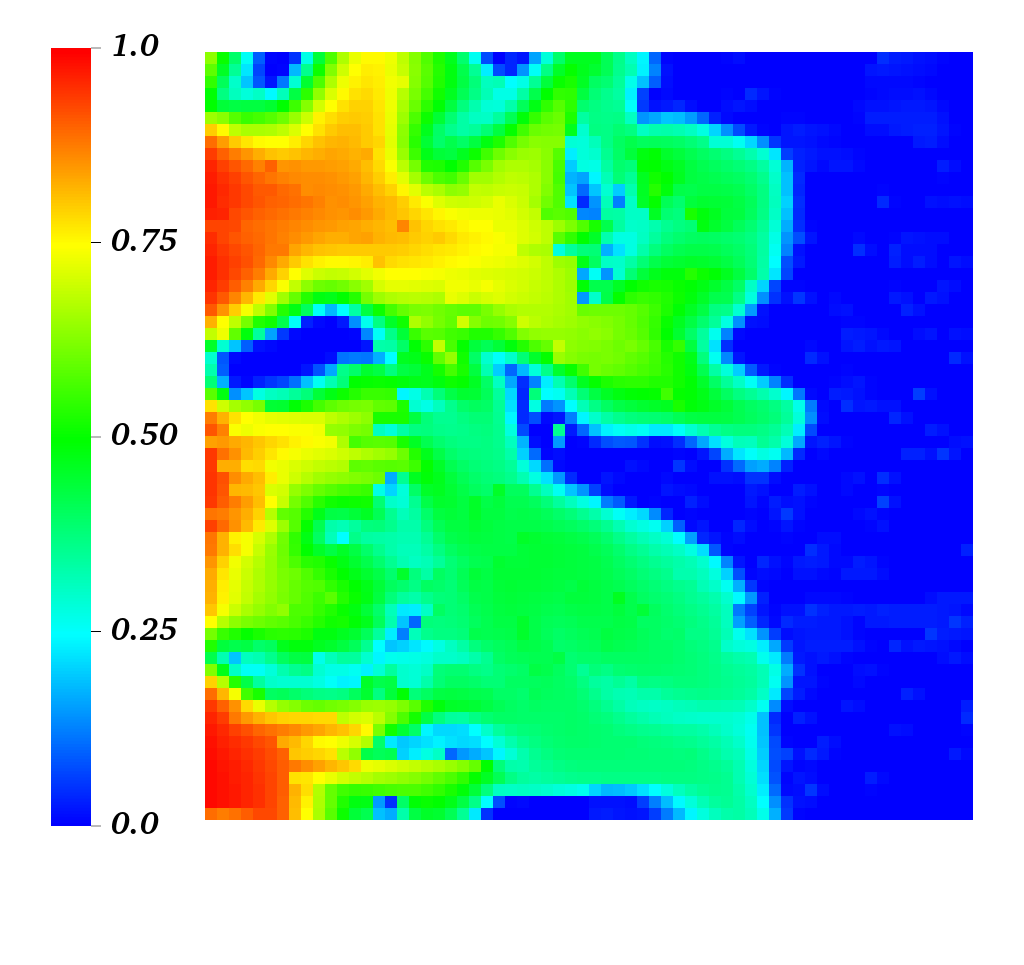}
        \caption{$10\%$ observations}
    \end{subfigure}
    \begin{subfigure}{0.3\textwidth}
        \includegraphics[width=\linewidth]{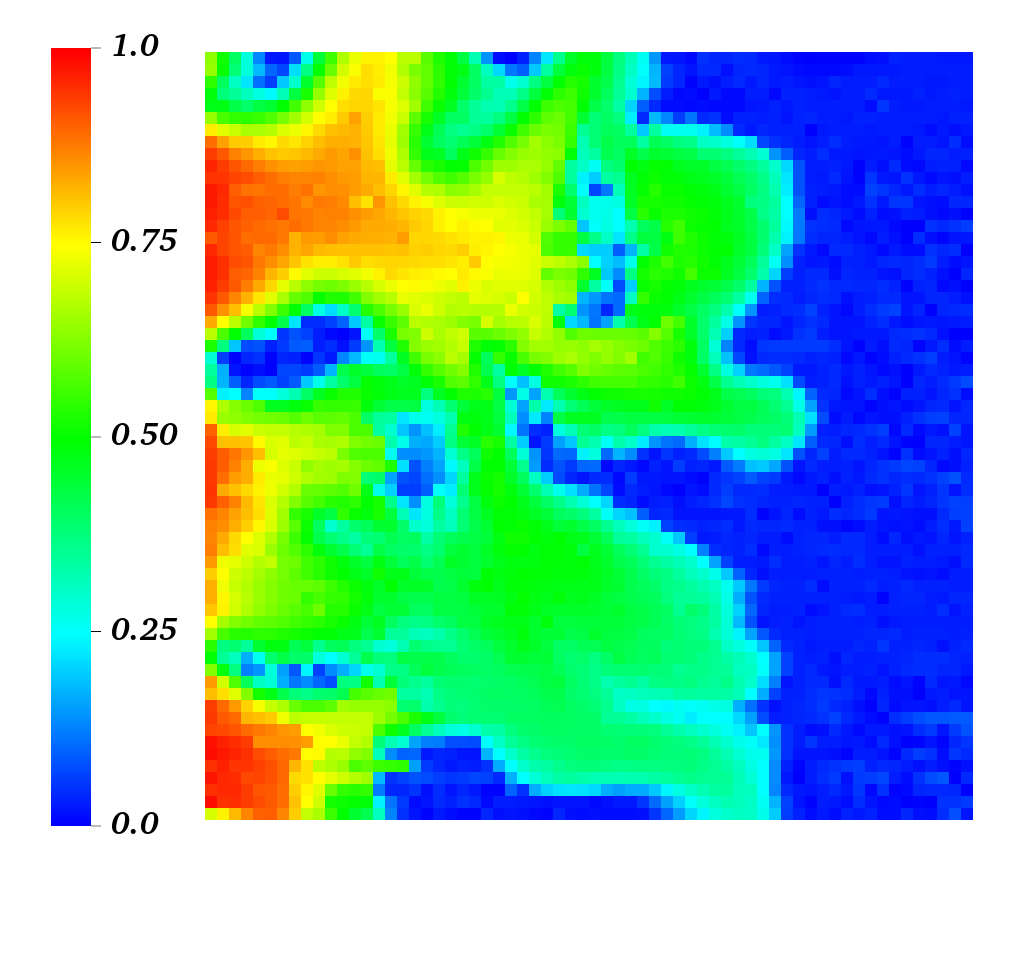}
        \caption{$30\%$ observations}
    \end{subfigure}
    \begin{subfigure}{0.3\textwidth}
        \includegraphics[width=\linewidth]{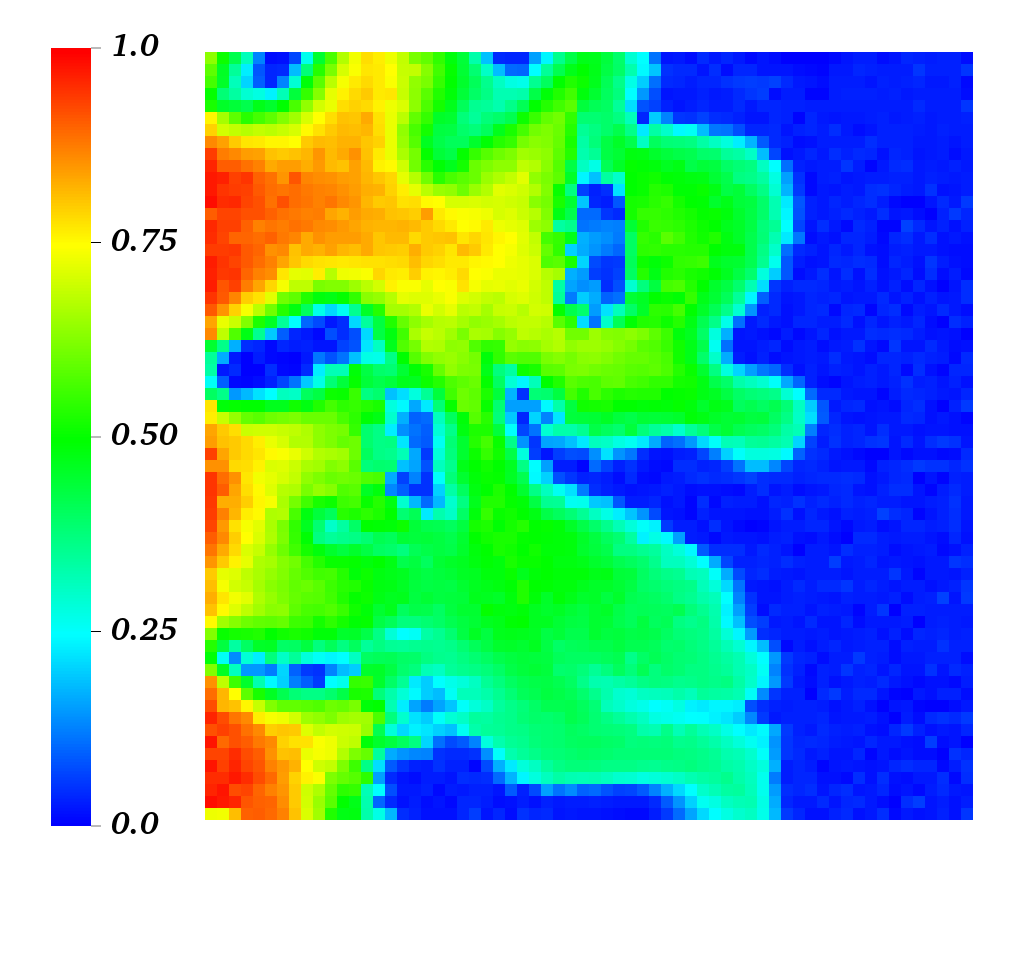}
        \caption{$50\%$ observations}
    \end{subfigure}
    \caption{EnSF estimated saturation ($\hat{s}_h$) with partial observations of the velocity, pressure, and saturation data.  }\label{Ex4: EnSF_S}
\end{figure}

\begin{figure}[h!]
    \centering
    \begin{subfigure}{0.3\textwidth}
        \includegraphics[width=\linewidth]{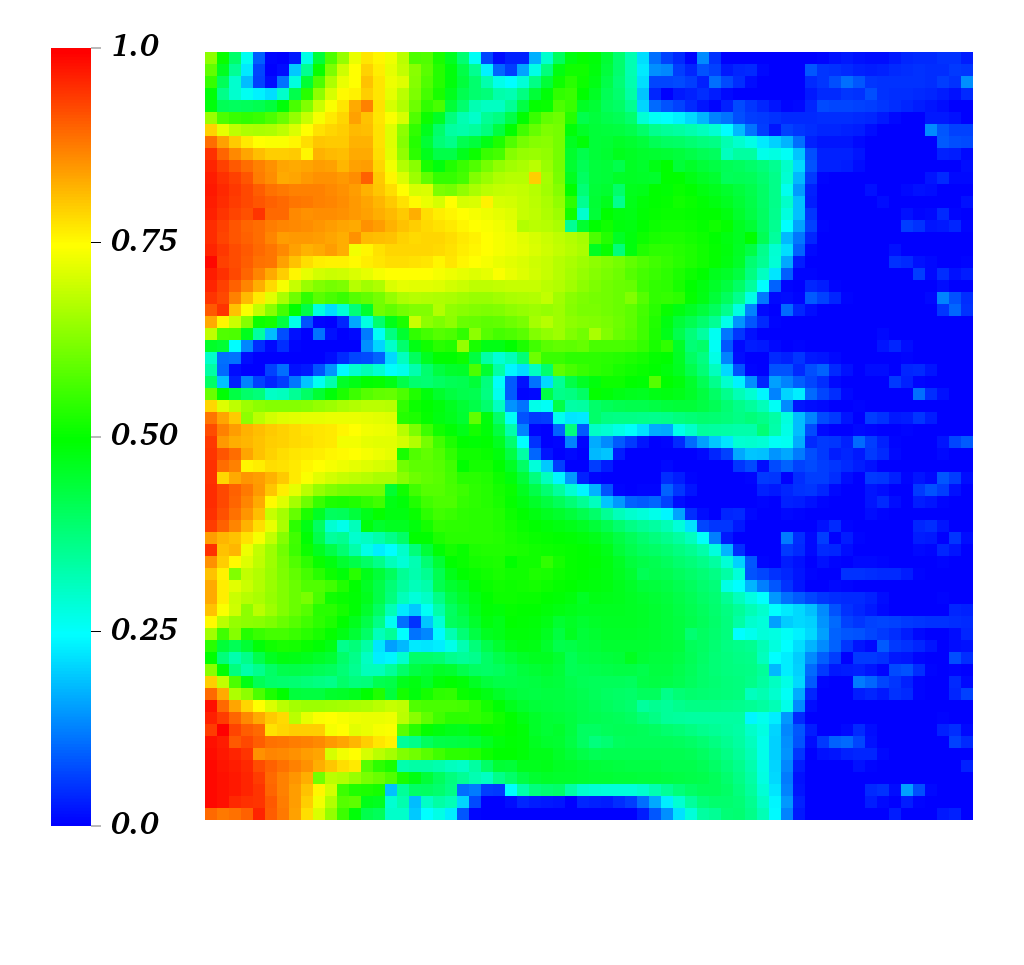}
        \caption{$10\%$ observations.}
    \end{subfigure}
    \begin{subfigure}{0.3\textwidth}
        \includegraphics[width=\linewidth]{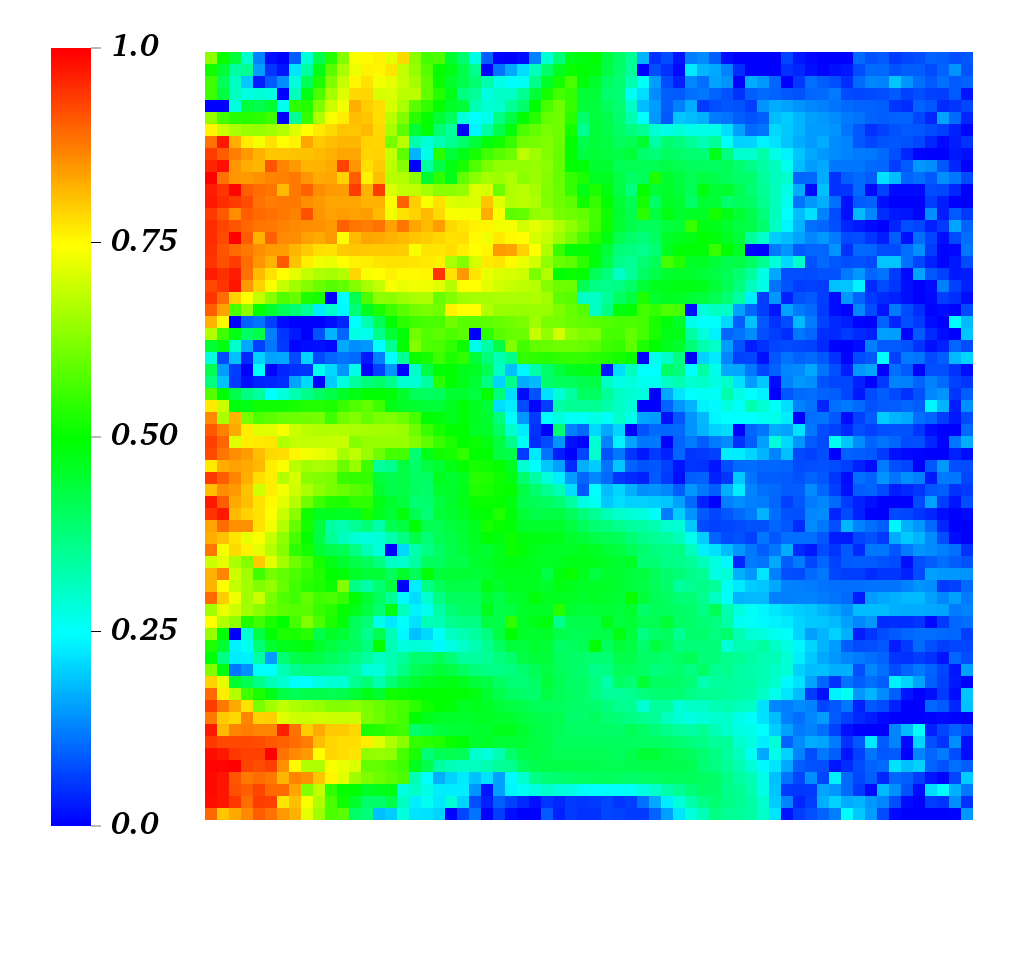}
        \caption{$30\%$ observations.}
    \end{subfigure}
    \begin{subfigure}{0.3\textwidth}
        \includegraphics[width=\linewidth]{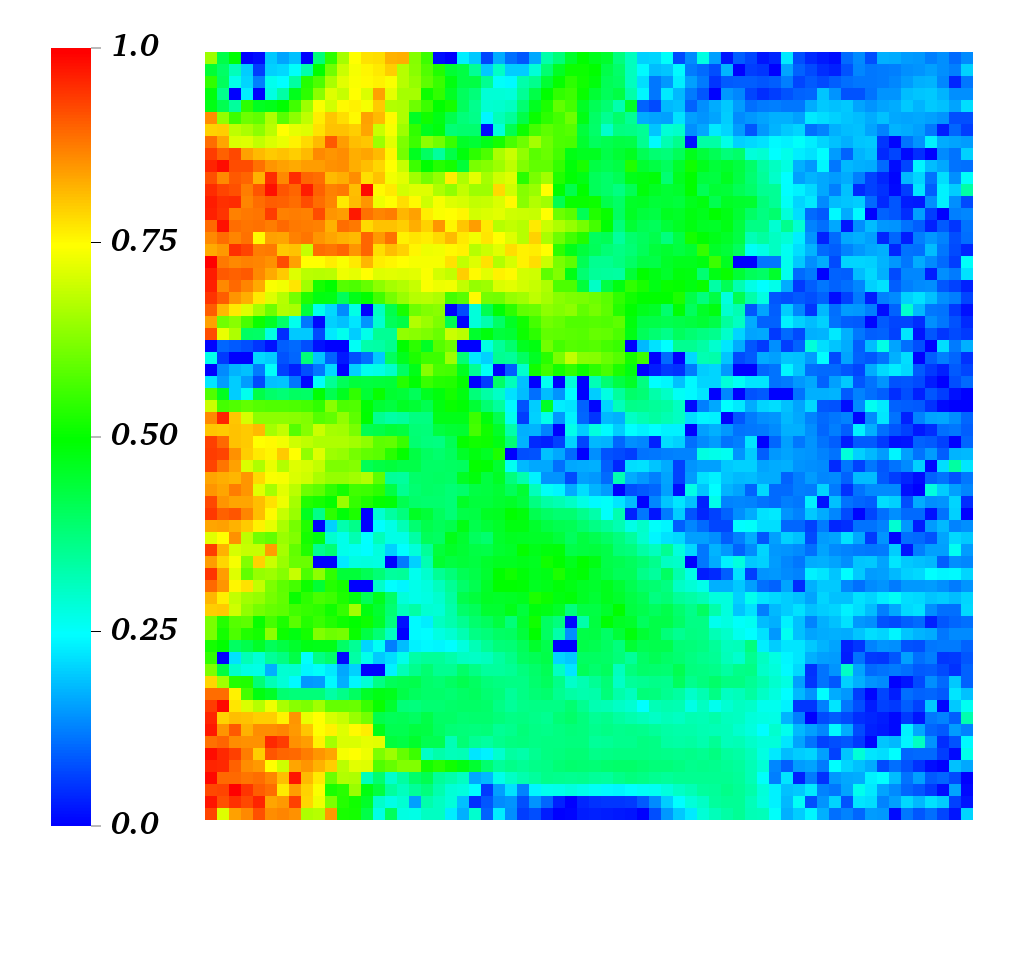}
        \caption{$50\%$ observations.}
    \end{subfigure}
    \caption{LETKF estimation of saturation ($\hat{s}_h$) with partial observations of the velocity, pressure, and saturation data. }\label{Ex4: LETKF_S}
\end{figure}

In Fig. \ref{Ex4: EnSF_S} and Fig. \ref{Ex4: LETKF_S}, we present the estimated saturation obtained using the EnSF and the LETKF, respectively. Each figure illustrates the results under three scenarios, where $10\%$, $30\%$, and $50\%$ of the state variables (velocity, pressure, and saturation) are observed. These figures demonstrate that the EnSF generally outperforms the LETKF in estimating the saturation. Notably, with only $10\%$ of observational data, the EnSF can still roughly identify the presence of low-permeability regions. As the proportion of observations increases to $30\%$ and $50\%$, the EnSF accurately captures the features of the reference solution for saturation, including the low-permeability regions. \par

In contrast, the LETKF provides less accurate saturation estimates compared to the EnSF. Notably, the LETKF struggles to identify the low-permeability regions, as it relies heavily on the accuracy of the prediction model \cite{AEnKF}. When the prediction model deviates significantly from the ground-truth data that reflect the actual state of the system, the LETKF's performance deteriorates. Additionally, the LETKF's limitations in handling nonlinear observational data further contribute to its poorer performance. Specifically, the inclusion of $20\%$ $\arctan$ observations introduces additional challenges, leading to noisier LETKF-based estimates.\par
\begin{figure}[h!]
    \centering
    \begin{subfigure}{0.3\textwidth}
        \includegraphics[width=\linewidth]{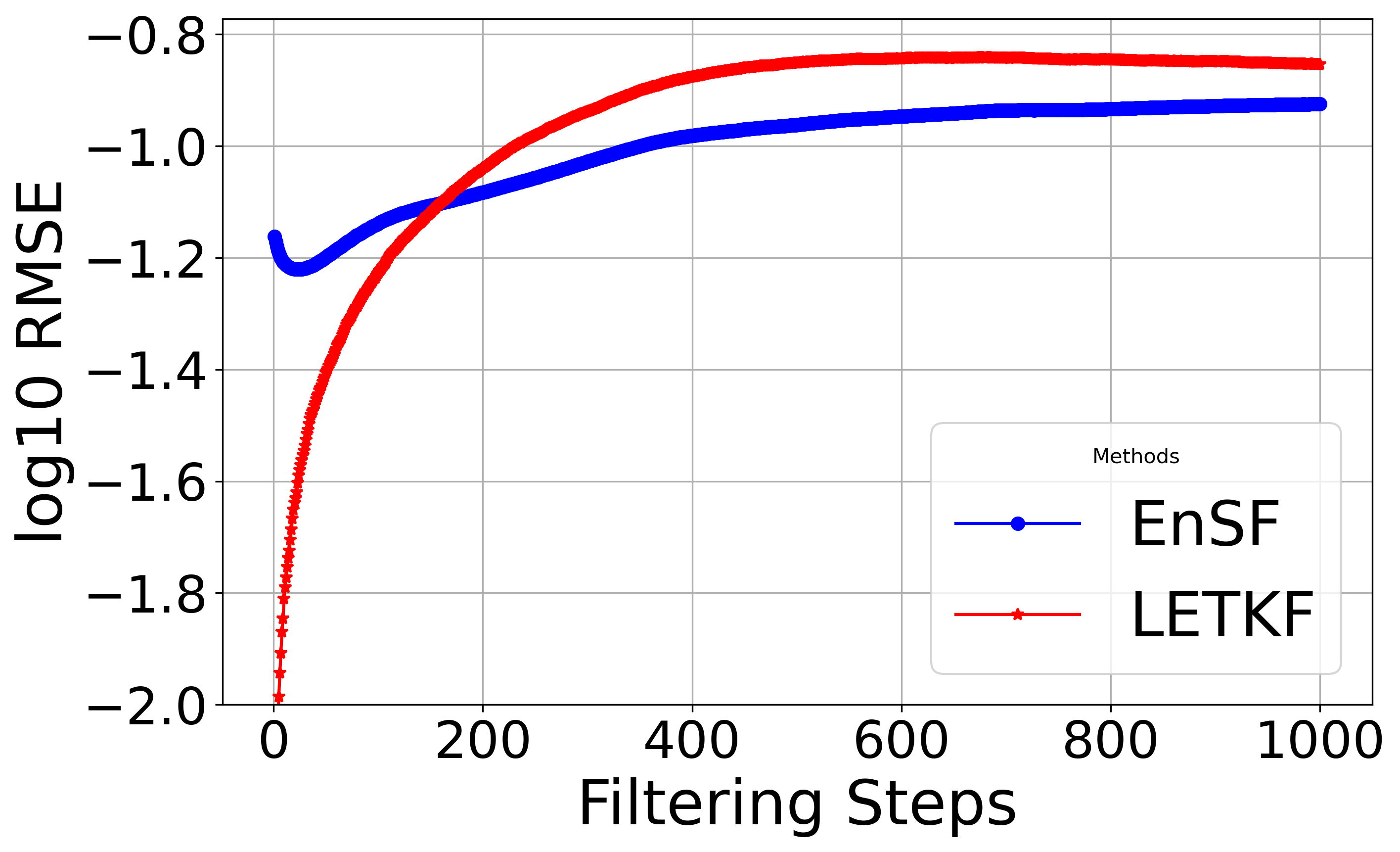}
        \caption{RMSEs comparison with 10\% observations}
        \label{Ex4:RMSEsa}        
    \end{subfigure}
    \begin{subfigure}{0.3\textwidth}
        \includegraphics[width=\linewidth]{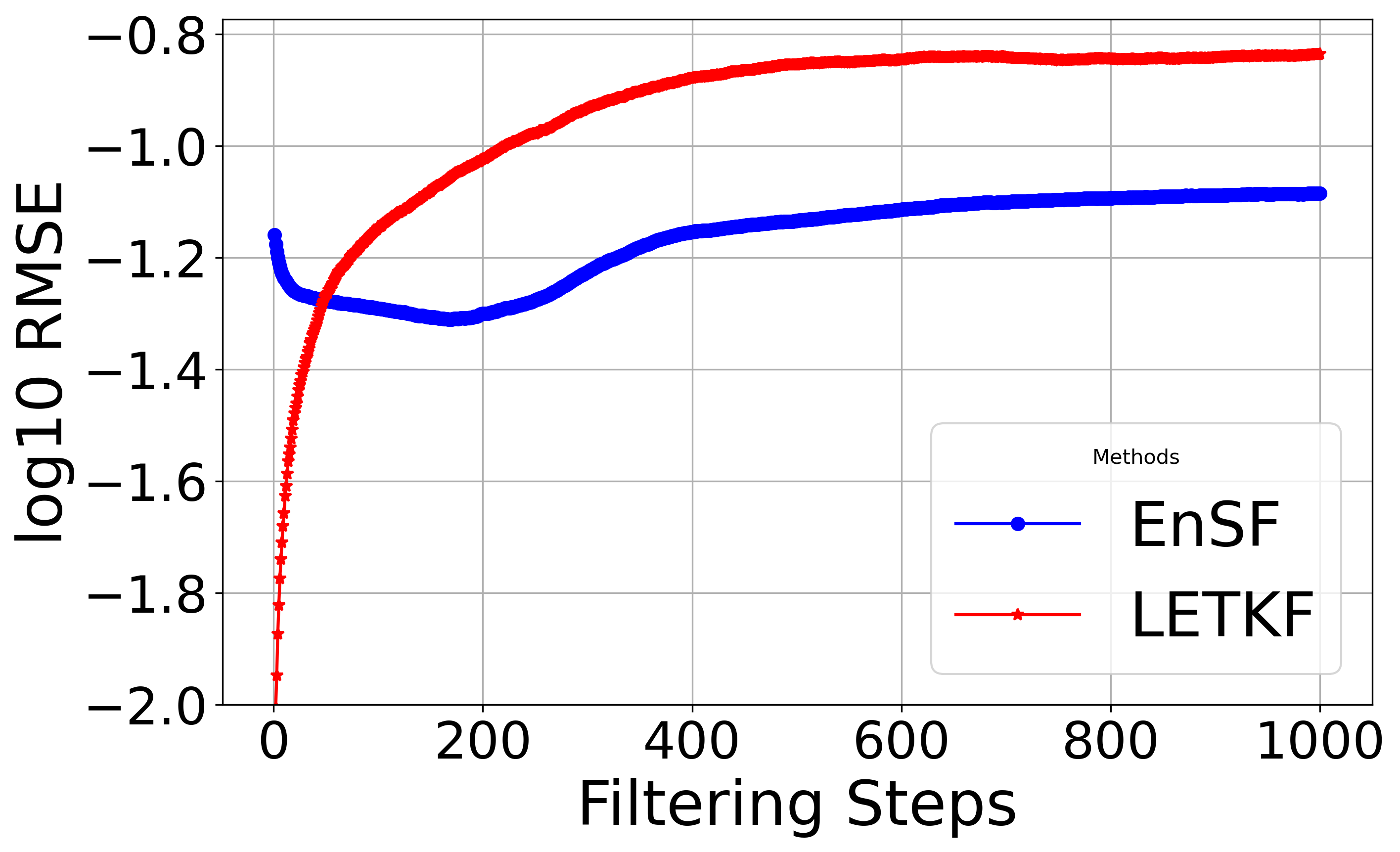}
        \caption{RMSEs comparison with 30\% observations}
        \label{Ex4:RMSEsb}        
    \end{subfigure}
    \begin{subfigure}{0.3\textwidth}
        \includegraphics[width=\linewidth]{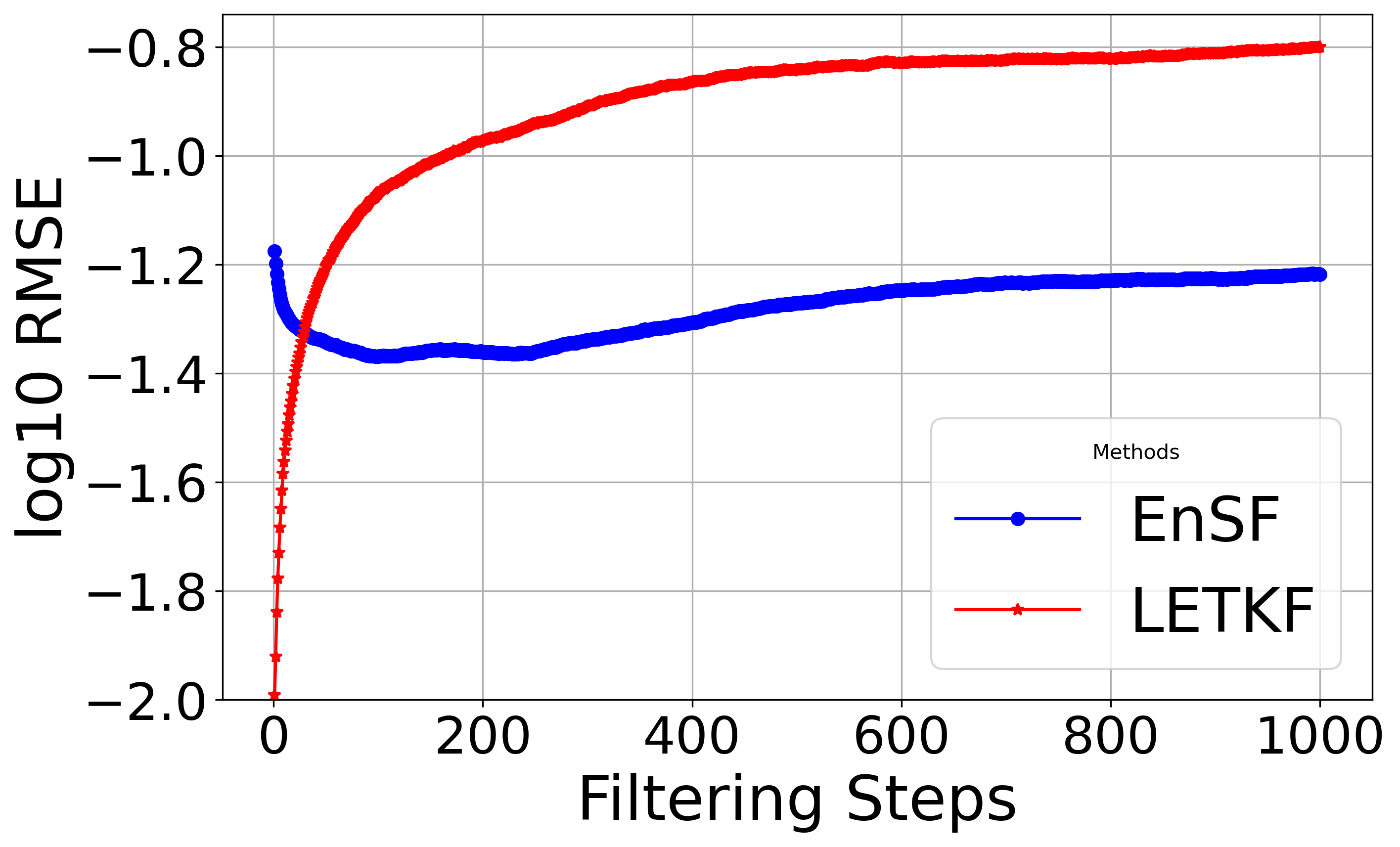}
        \caption{RMSEs comparison with 50\% observations}
        \label{Ex4:RMSEsc}        
    \end{subfigure}
    \caption{Comparison of RMSEs with $10\%$, $30\%$, and $50\%$ state observations. With only $10\%$ of observational data, the two methods have similar performance. With more and more data available, the EnSF starts to outperform the LETKF significantly. }
    \label{Ex4:RMSEs}
\end{figure}

To construct a more straightforward comparison, we calculate the RMSEs for each method by aggregating the estimation errors across all mesh grid points for $\tilde{X}$. We then plot the average RMSEs over time steps. The RMSEs comparisons for 10\%, 30\%, and 50\% state observations are presented in Fig. \ref{Ex4:RMSEsa}, \ref{Ex4:RMSEsb}, and \ref{Ex4:RMSEsc}, respectively. From these figures, we observe that with only 10\% of observational data, the EnSF and the LETKF show similar accuracy, as limited data makes it difficult to effectively calibrate the state estimation. However, as more data becomes available, the EnSF starts to outperform the LETKF significantly. This example highlights the strength of our EnSF-based data assimilation for capturing low-permeability regions compared to the LETKF method. 

\begin{figure}[ht]
    \centering
    \includegraphics[width=0.4\linewidth]{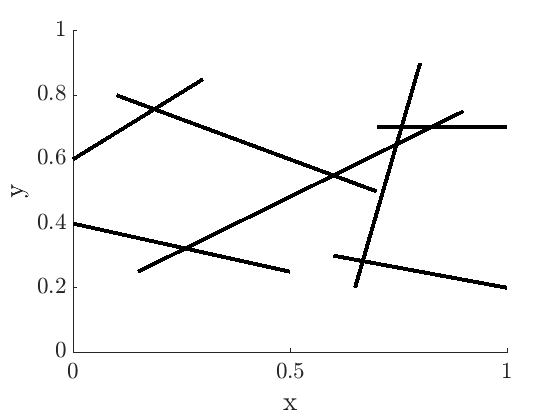}
    \caption{Complete Fracture Network}
    \label{fig:fracture}
\end{figure}
To generate the reference data, we define the following $k(\bx)$:
\begin{equation}
    k(\bx) := \min\left\{ \exp(-400 d_s^2) + \eta_1(\bx) + \eta_2(\bx) + \eta_3(\bx) + \eta_4(\bx), 4 \right \},
\end{equation}
where $d_s$ is shortest perpendicular distance of a point $\bx$ from the fracture. $\eta_i(\bx)\}_{i=1}^3$ represent the uncertainties defined as follows:
\begin{itemize}
\item $\eta_1$ occurs in $7\%$ of the region, following $\eta_1(\bx) \sim N(1.5, 0.75)$, 
\item $\eta_2$ occurs in $2\%$ of the region, following $\eta_2(\bx) \sim N(1.0, 0.5)$,
\item $\eta_3$ occurs in $1\%$ of the region, following $\eta_3(\bx) \sim N(2.0, 1.0)$.
\end{itemize}
\subsection{Exampe 3. A case with a fracture network}
In this example, we consider a scenario where a fracture network exists within a domain, but we have limited information about the fracture, i.e., its precise location and structure. These fractures represent the region of high permeability, with high saturation flow along them. With the limited information the fractures, the simulation in our model deviates largely from the ground truth. Using the EnSF method, we show that even with incomplete knowledge of the fracture network, we can calibrate the saturation from the partial and noisy observations of the data. The corrections not only improve our predictions but also provide valuable insights into the fracture network. This mirrors a real-world scenario where the exact underground fracture network is unknown. However, we may have access to sensors that collect the saturation data from various locations. \par
We assume that the actual fracture is as shown in Fig. \ref{fig:fracture}.

\begin{figure}[h!]
    \centering
    \begin{subfigure}{0.3\textwidth}
        \includegraphics[width=\linewidth]{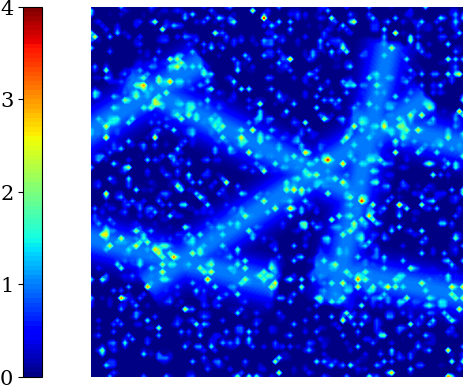}
        \caption{$k(\bx)$}
        \label{fig:Ex3_kx}        
    \end{subfigure}
    \begin{subfigure}{0.3\textwidth}
        \includegraphics[width=\linewidth]{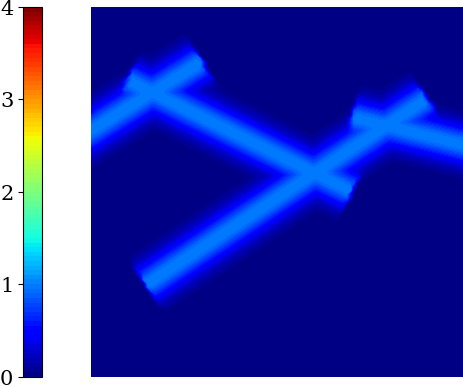}
        \caption{$\hat{k}(\bx)$}
        \label{fig:Ex3_k_hat_x}        
    \end{subfigure}
    \caption{Permeability function }
    \label{fig:Ex3_k}
\end{figure}

With the above definition, the permeability function is shown in Fig. \ref{fig:Ex3_kx}. For the forward model $f(\cdot)$, we use the permeability $\hat{k}(\bx)$ that only has partial information of the fracture network. Fig. \ref{fig:Ex3_k_hat_x} shows the permeability function for the forward model.

\begin{figure}[h!]
    \centering
    \begin{subfigure}[t]{0.3\textwidth}
        \includegraphics[width=\linewidth]{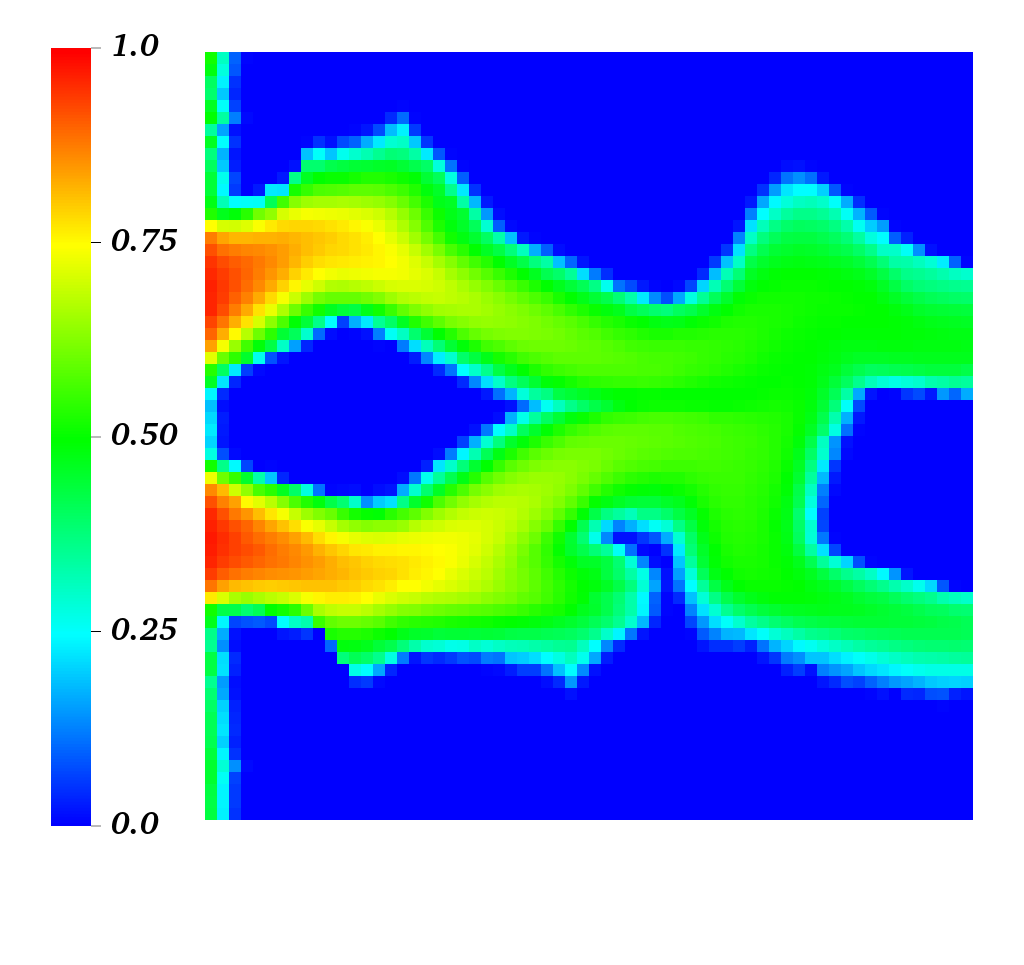}
        \caption{Reference saturation $(\hat{s}_h)$ using $k(\bx)$}
        \label{Ex5a: ref}
    \end{subfigure}
    \begin{subfigure}[t]{0.3\textwidth}
        \includegraphics[width=\linewidth]{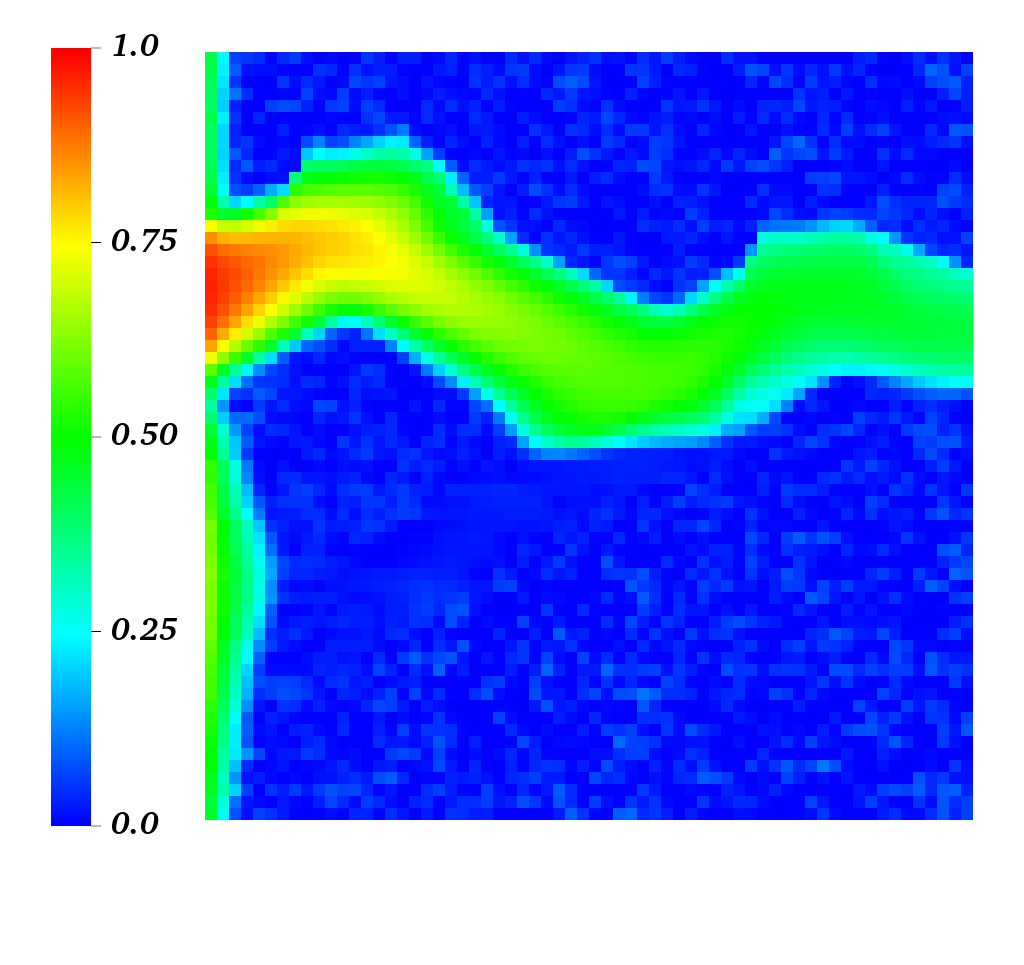}
        \caption{Estimated saturation with no observation using $\hat{k}(\bx)$}
        \label{Ex5b: ref}
    \end{subfigure}    
    \begin{subfigure}[t]{0.3\textwidth}
        \includegraphics[width=\linewidth]{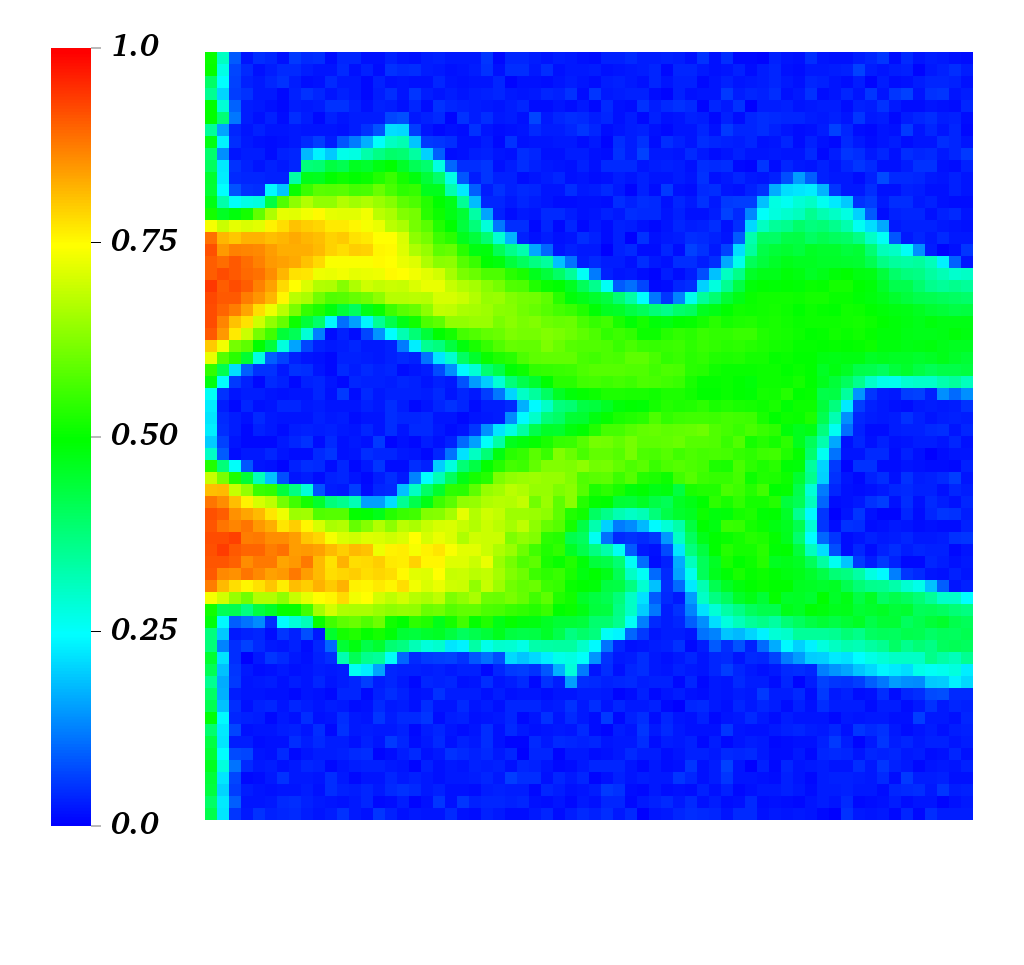}
        \caption{$100\%$ observations by using EnSF}
        \label{Ex5c: ref}
    \end{subfigure}
    \caption{The comparison between reference saturation, estimated saturation with full data assimilation, and estimated saturation without data assimilation at $t = 1$.  }\label{Ex5: ref}
\end{figure}

\begin{figure}[h!]
    \centering
    \begin{subfigure}[t]{0.3\textwidth}
        \includegraphics[width=\linewidth]{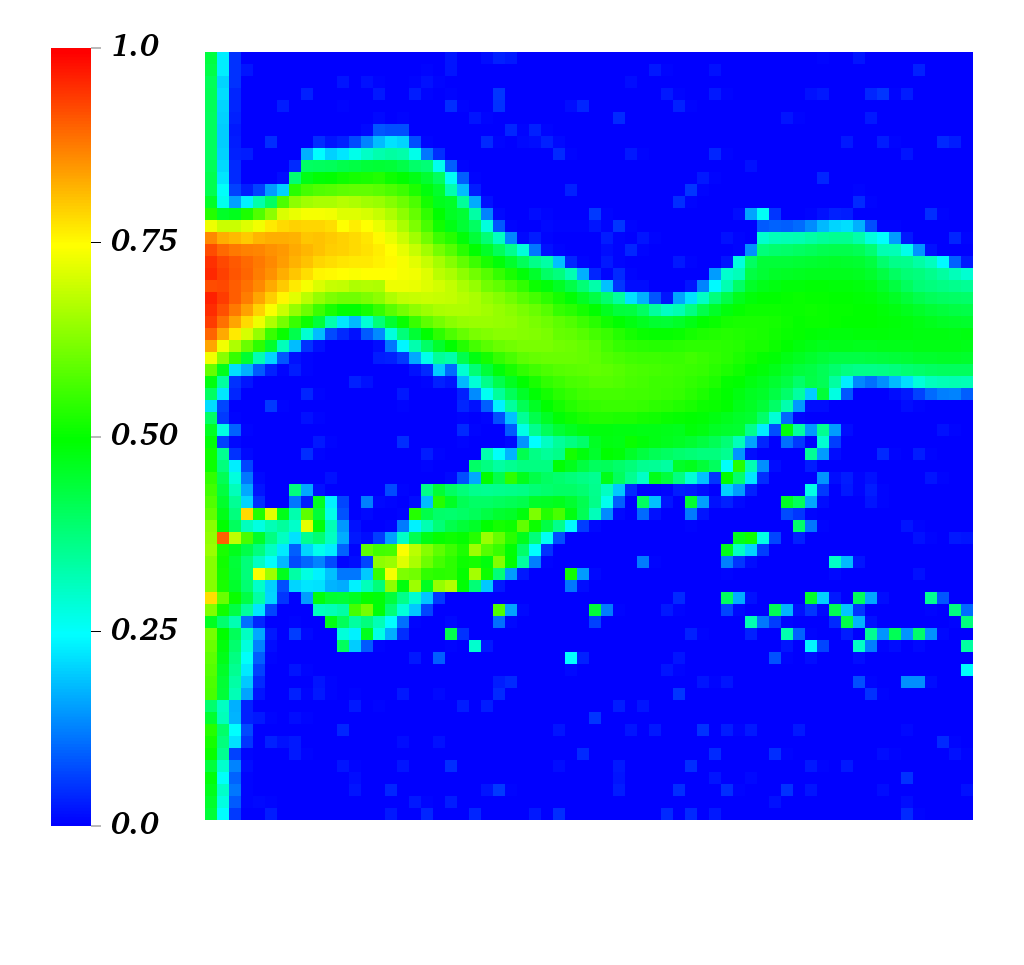}
        \caption{$10\%$ observations by EnSF.}
    \end{subfigure}
    \begin{subfigure}[t]{0.3\textwidth}
        \includegraphics[width=\linewidth]{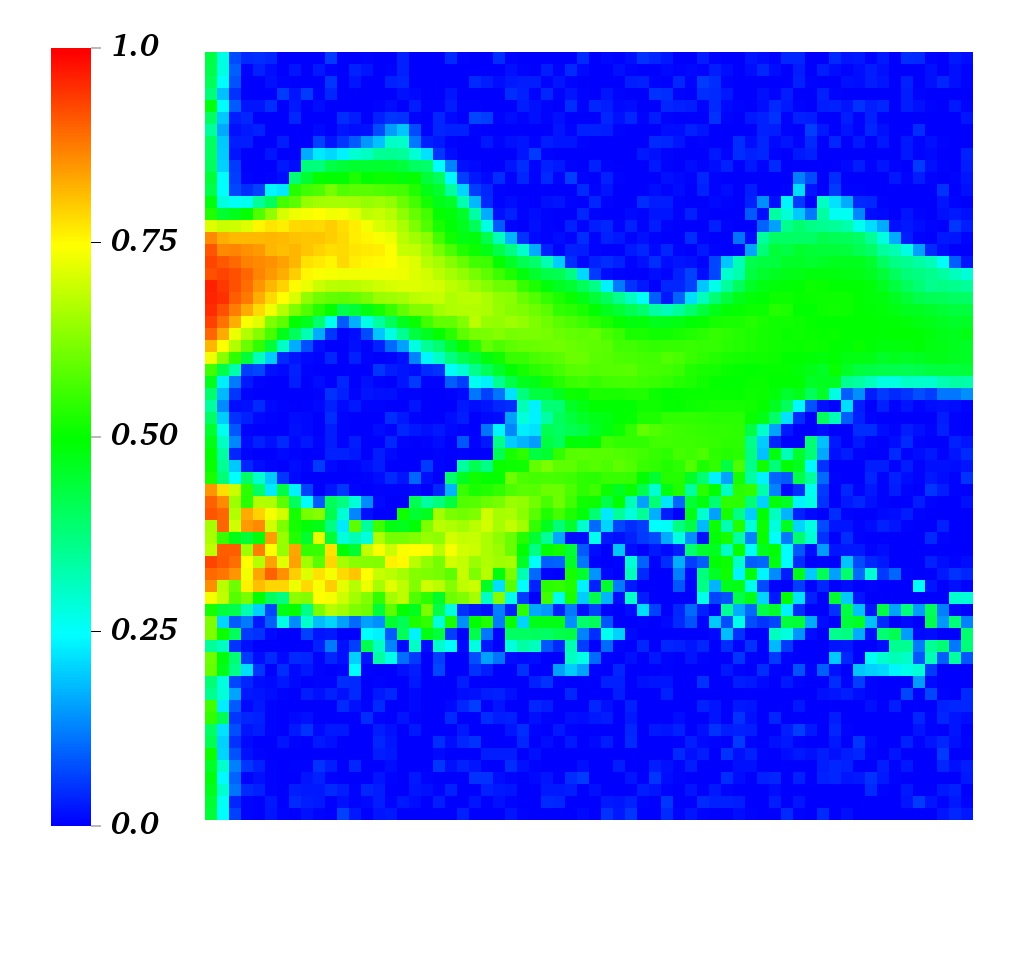}
        \caption{$40\%$ observations by EnSF.}
    \end{subfigure}
    \begin{subfigure}[t]{0.3\textwidth}
        \includegraphics[width=\linewidth]{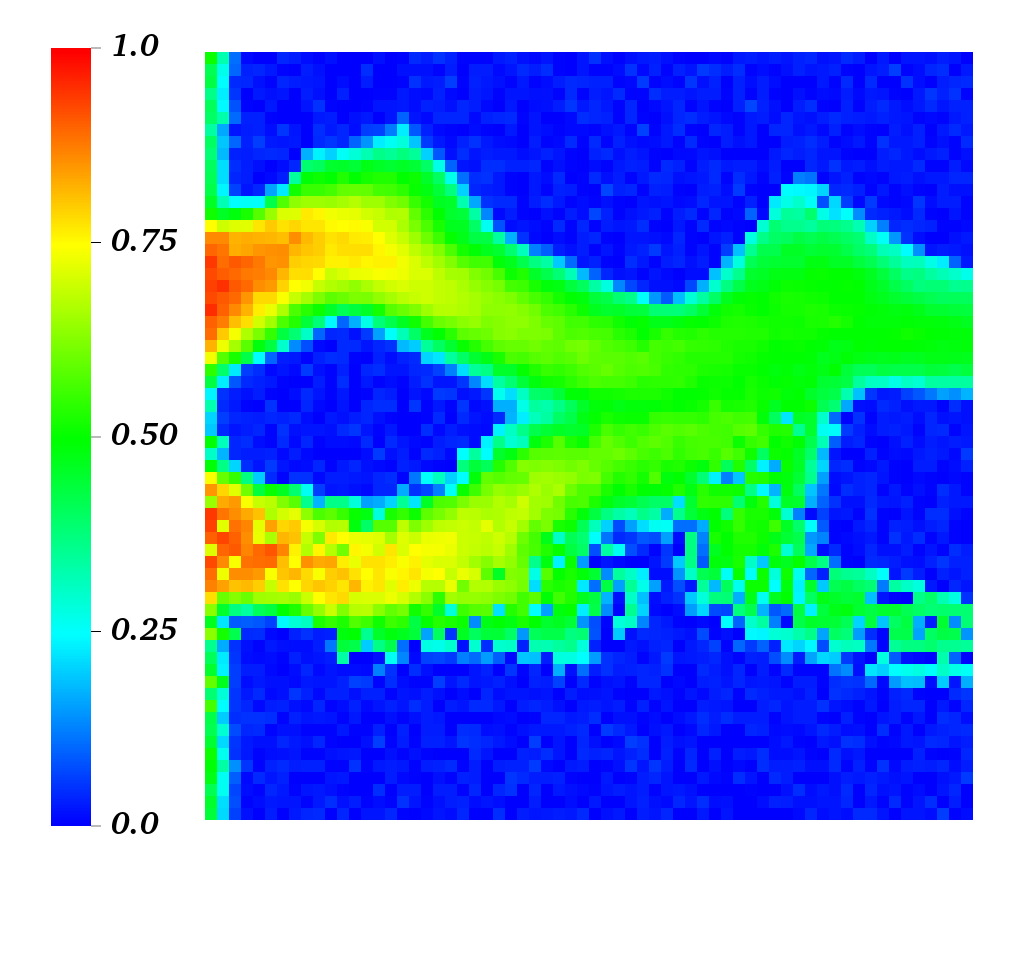}
        \caption{$70\%$ observations by EnSF.}
    \end{subfigure}
    \begin{subfigure}[t]{0.3\textwidth}
        \includegraphics[width=\linewidth]{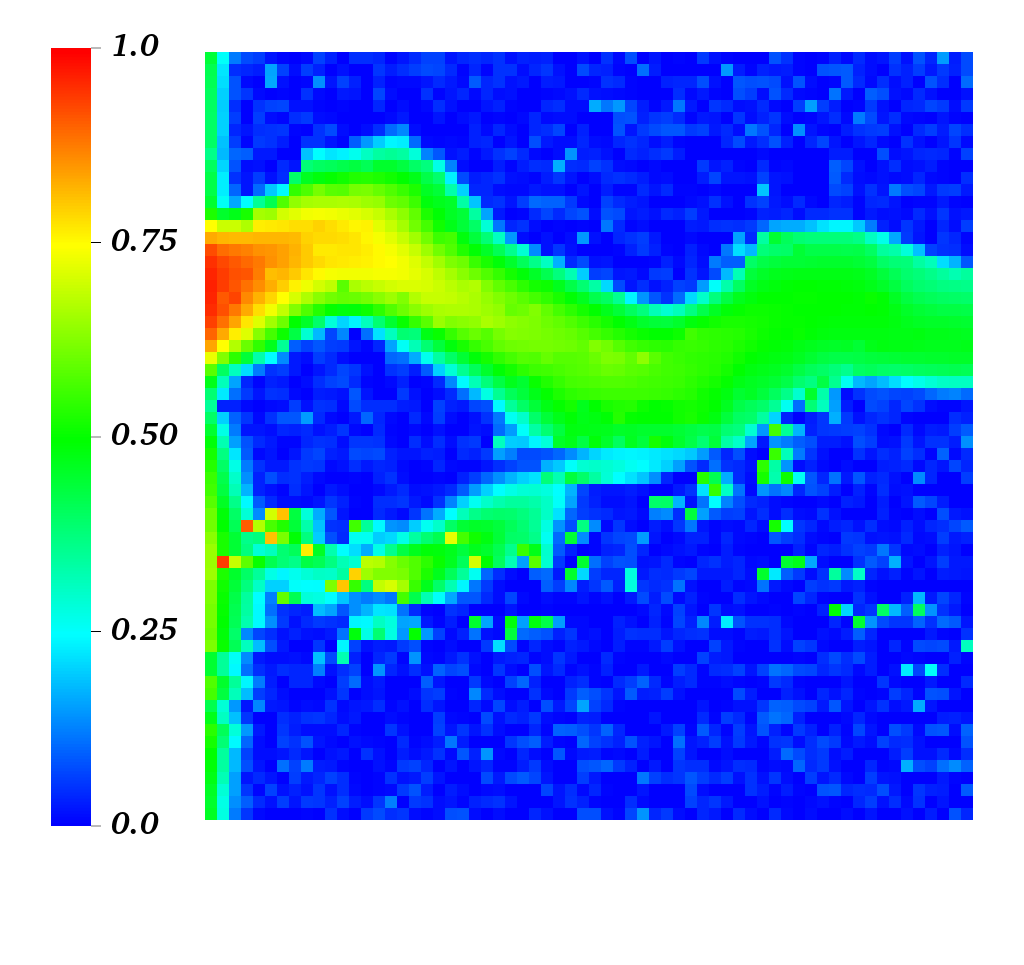}
        \caption{$10\%$ observations by LETKF.}
    \end{subfigure}
    \begin{subfigure}[t]{0.3\textwidth}
        \includegraphics[width=\linewidth]{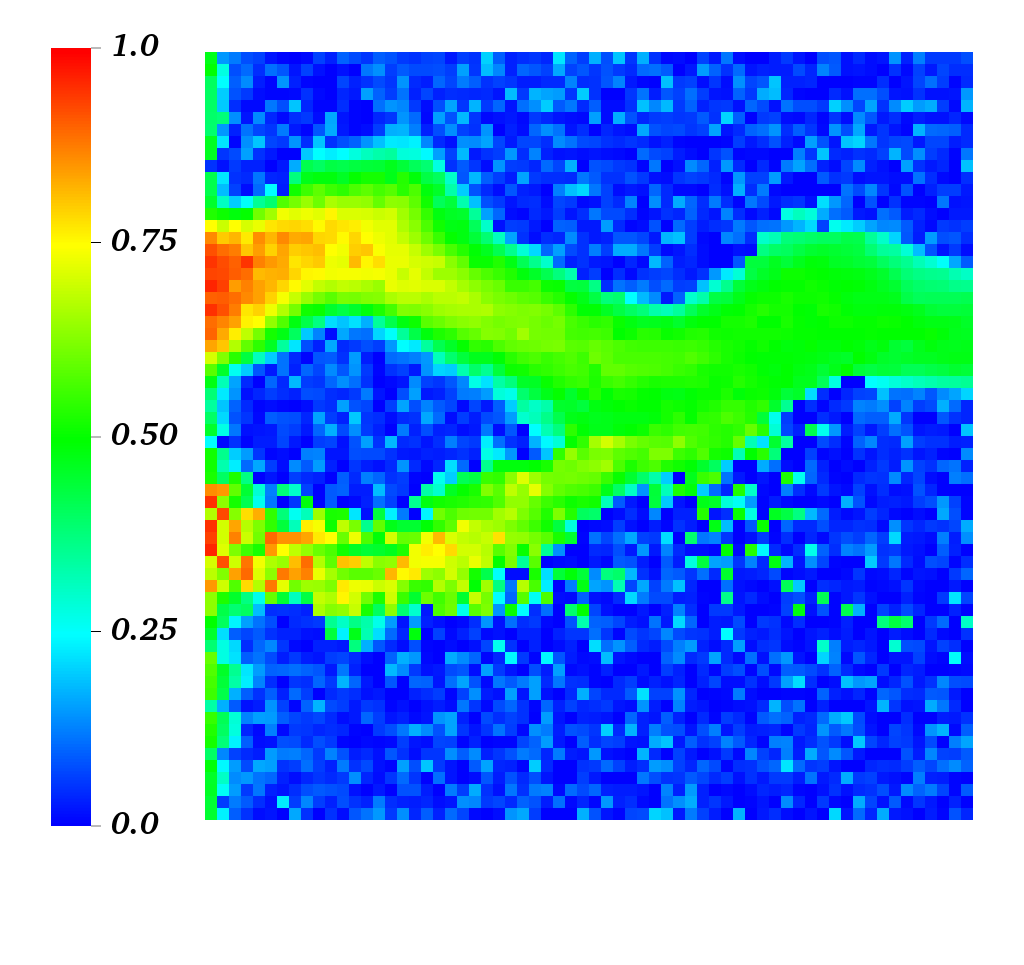}
        \caption{$40\%$ observations by LETKF.}
    \end{subfigure}
    \begin{subfigure}[t]{0.3\textwidth}
        \includegraphics[width=\linewidth]{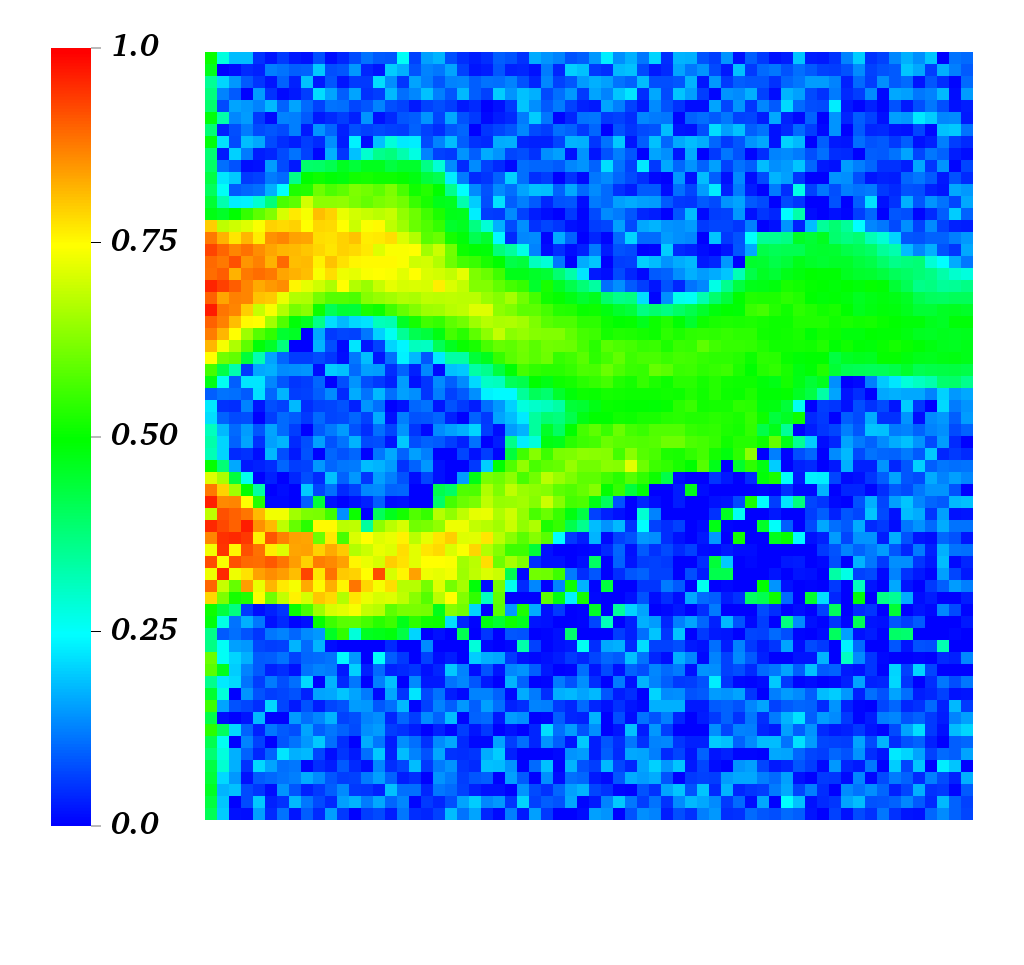}
        \caption{$70\%$ observations by LETKF.}
    \end{subfigure}
    \caption{EnSF and LETKF estimated saturation $(\hat{s}_h)$ with partial observations of the velocity, pressure and saturation at $t = 1$.  }\label{Ex5: partial}
\end{figure}

\begin{figure}[h!]
    \centering
    \begin{subfigure}{0.3\textwidth}
        \includegraphics[width=\linewidth]{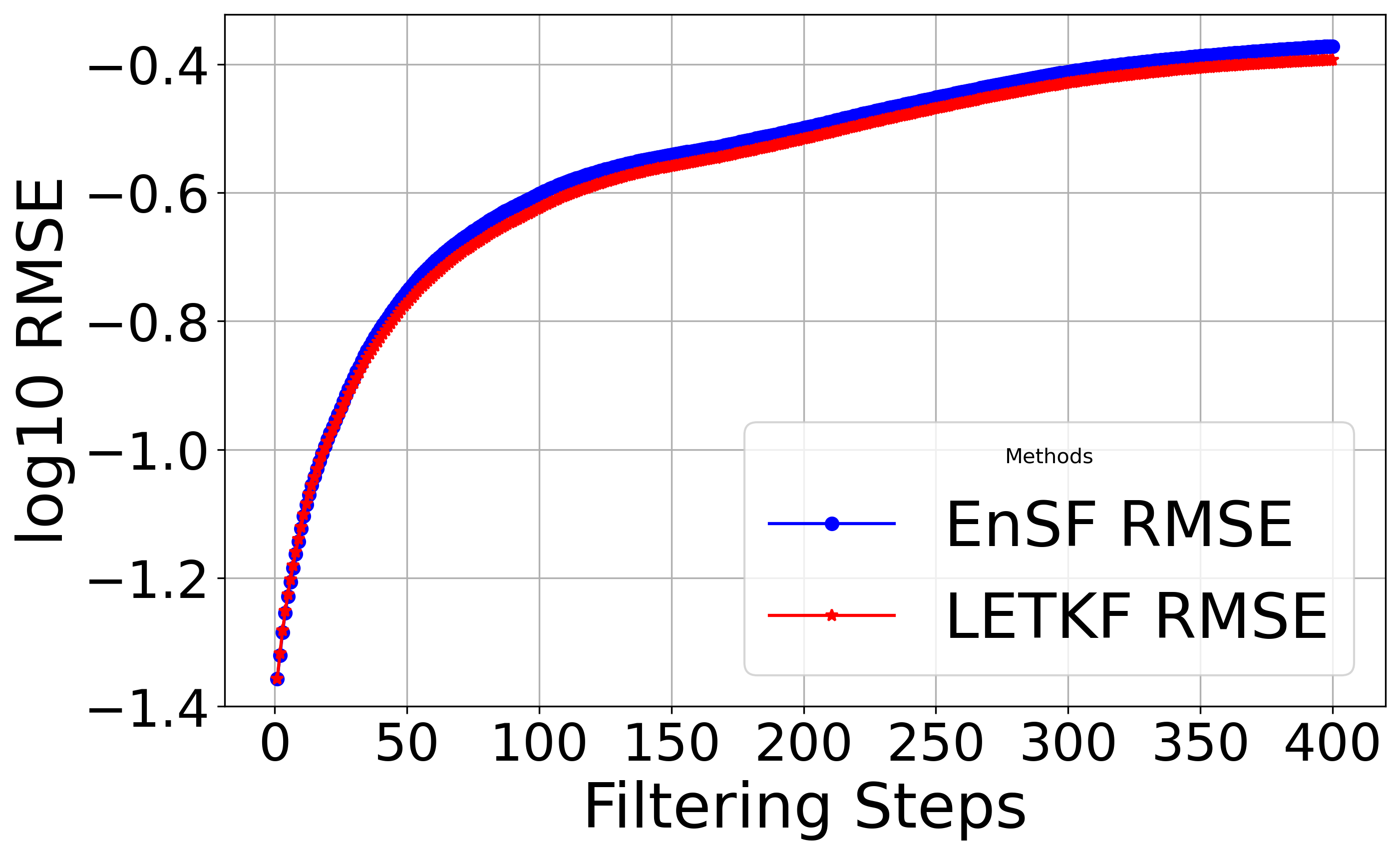}
        \caption{RMSEs comparison with 10\% observations}
        \label{Ex5:RMSEsa}        
    \end{subfigure}
    \begin{subfigure}{0.3\textwidth}
        \includegraphics[width=\linewidth]{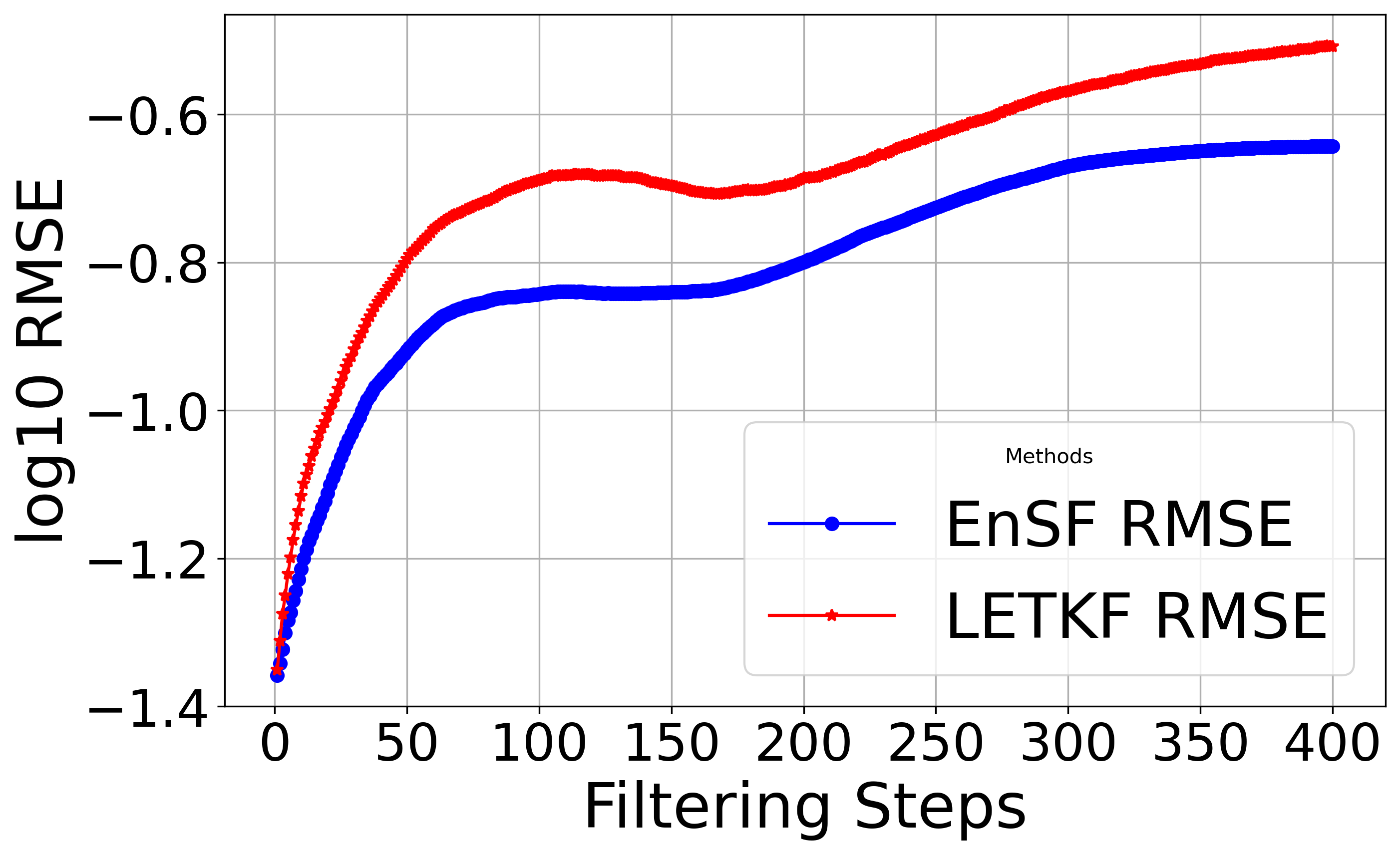}
        \caption{RMSEs comparison with 40\% observations}
        \label{Ex5:RMSEsb}        
    \end{subfigure}
    \begin{subfigure}{0.3\textwidth}
        \includegraphics[width=\linewidth]{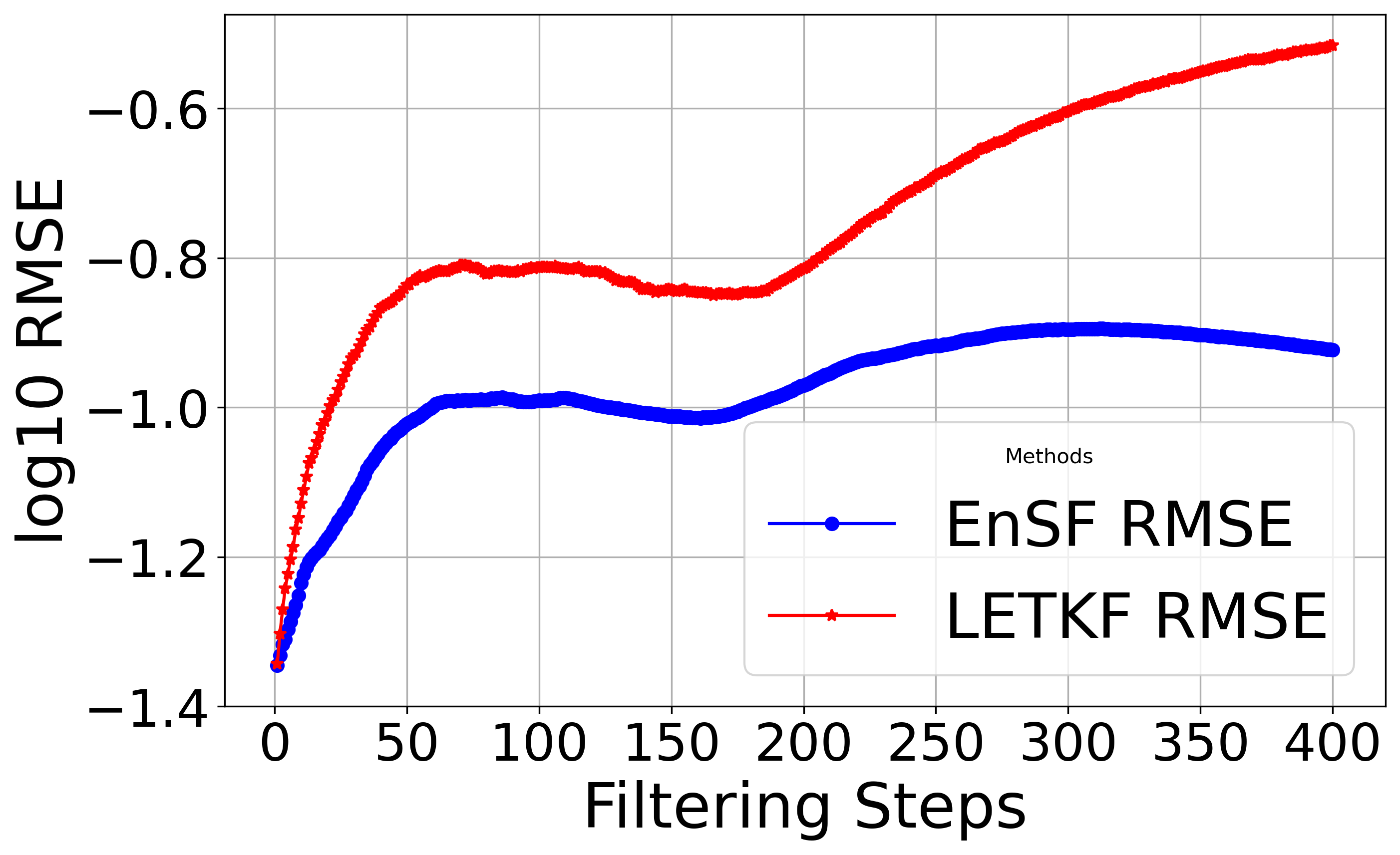}
        \caption{RMSEs comparison with 70\% observations}
        \label{Ex5:RMSEsc}        
    \end{subfigure}
    \caption{Comparison of RMSEs with $10\%$, $40\%$, and $70\%$ state observations. With only $10\%$ of observational data, the two methods have similar performance. With more and more data available, the EnSF starts to outperform the LETKF significantly. }
    \label{Ex5:RMSE}
\end{figure}

Similar to the previous example, we fix the randomness, $\{\eta_i(\bx)\}_{i=1}^3$, in permeability, and we perturb the initial condition. The data assimilation procedure in this example is implemented over the time interval $(0, 1]$ with temporal step-size $\delta t = 0.0025$ and observation defined in Eqn.\eqref{NL}. \par
In Fig. \ref{Ex5: ref}, we compare the results with and without our data assimilation algorithm. Fig. \ref{Ex5a: ref} shows the reference data, saturation $\hat{s}_h$, at $t = 1$. Fig. \ref{Ex5b: ref} is the final saturation obtained without the data assimilation, i.e., the saturation at $t = 1$ by the state dynamics $\bar{f}(\cdot)$ using $\hat{k}(\bx)$ with a random initial condition. Since the forward model lacks information on some major branches of the fracture, we observe that the flow of saturation is only through the fracture defined in $\hat{k}(\bx)$. On the other hand, using the EnSF-based data assimilation with the full observation of the saturation, we recover the flow through the missing fractures. Fig. \ref{Ex5c: ref} shows that the result is in agreement with the reference solution. This shows that with the data assimilation method, we can calibrate the saturation with limited knowledge on fractures. \par

In Fig. \ref{Ex5: partial}, we plot the EnSF and LETKF estimated saturation under $10\%, 40\%, 70\%$ observations. 
%$10\%, 40\%, 70\%$ observations on areas with significant different low-permeability area as shown in Fig. \ref{fig:Ex3_k}, that is the regions that without fractures in the forward model but with fractures in the reference model. 
Similarly to what we observe in Example 2, Fig. \ref{Ex5: partial} shows that the EnSF has a better performance over LETKF in this task, which is more significant when the observation percentage increases. Fig. \ref{Ex5:RMSE} shows the RMSEs calculated by EnSF and LETKF in the region of fractures where the forward model has no former knowledge. It is clear that EnSF outperforms LETKF in calculating the solution located within fractures that are only observable in the observations, particularly as the observation percentage goes up from $10\%$ to $70\%$. Both Fig. \ref{Ex5: partial} and Fig. \ref{Ex5:RMSE} show the competency of EnSF in helping calibrate solutions in complicated underground fracture scenarios.

%\clearpage
%\section{Conclusions}

%\section*{Acknowledgments}
%The authors R. Hu, F. Bao, and S. Lee are partially supported by U.S. Department of Energy, Office of Science, Energy Earthshots Initiatives under Award DE-SC-0024703. F. Bao would also like to acknowledge the support from U.S. National Science Foundation through project DMS-2142672 and the support from the U.S. Department of Energy, Office of Science, Office of Advanced Scientific Computing Research, Applied Mathematics program under Grants DE-SC0025412.
%
\bibliographystyle{elsarticle-num} 

%\bibliography{Reference}

\end{document}